\documentclass[twocolumn]{svjour3mod}          
\smartqed  
\journalname{Computational Mechanics}

\usepackage[dvipsnames]{xcolor}
\usepackage{graphicx}
\usepackage{bigints}
\usepackage{amssymb}
\usepackage{amsmath}
\usepackage{bm}
\usepackage{mathtools, cuted}
\usepackage{subfigure}
\usepackage[hidelinks]{hyperref}

\usepackage{tikz}
\usepackage{pgfplots}
\usetikzlibrary{matrix}
\usepgfplotslibrary{groupplots}
\usepackage{tikz-3dplot}
\usetikzlibrary{calc,intersections,shadows,patterns,decorations.pathmorphing,shapes,arrows,decorations.pathmorphing,external}
\tikzset{>=latex}

\newcommand{\mainfigurespath}{figures/}

\graphicspath{{\mainfigurespath}}
\DeclareGraphicsExtensions{.pdf,.png,.jpg}

\usepackage{adjustbox}
\usepackage[colorinlistoftodos, obeyFinal]{todonotes}

\DeclareFontFamily{OT1}{pzc}{}
\DeclareFontShape{OT1}{pzc}{m}{it}{<-> s * [1.10] pzcmi7t}{}
\DeclareMathAlphabet{\mathpzc}{OT1}{pzc}{m}{it}

\usepackage[misc]{ifsym}

\newcommand{\diff}{\text{d}}

\newcommand{\dOp}{\: \text{d}\hat{\Omega}}

\newcommand{\bx}{\bm{x}}

\newcommand\Matrix[1]{\mbox{\boldmath$\sf #1$}}

\pgfplotsset{basic_paper_plots_style/.style={%
   width=0.95\textwidth, height=0.59\textwidth,
   legend columns=3,
   legend style={anchor=south, at={(0.5, 1.05)}, draw=none, font=\small}}}

\definecolor{green2}{RGB}{154, 205, 50}
\definecolor{green3}{RGB}{141, 182, 0}
\definecolor{blue2}{RGB}{0, 102, 255}
\definecolor{lightgreen}{RGB}{178,255,102}
\definecolor{papergreen}{rgb}{0.0, 0.5, 0.0}

\newif\ifdraftversion
\draftversionfalse

\newcommand{\ic}{{\rm\kern.24em                  
   \vrule width.02em height0.9ex depth-.05ex
   \kern-.26em c}}
\newcommand{\IC}{{\rm\kern.24em                  
   \vrule width.02em height1.4ex depth-.05ex
   \kern-.26em C}}

\newcommand{\eps}{\varepsilon}

\pgfplotsset{basic line style/.style={
             mark options=solid, thick}}

\pgfplotsset{solid 1/.style={basic line style,
             blue}}

\pgfplotsset{solid 2/.style={basic line style,
             black, dashed,mark=x}}

\pgfplotsset{solid 3/.style={basic line style,
             black,mark=x}}

\pgfplotsset{solid 2 C0/.style={basic line style,
             black, dashed, mark=o}}

\pgfplotsset{ss ANS-like/.style={basic line style,
             blue,mark=triangle}}

\pgfplotsset{ss/.style={basic line style,
             red,mark=square}}

\pgfplotsset{ss C0/.style={basic line style,
             red, mark=o}}

\pgfplotsset{ss red/.style={basic line style,
             violet, mark=x}}

\pgfplotsset{ANS/.style={basic line style,
             papergreen, mark=star}}

\pgfplotsset{reference/.style={basic line style,
             black, dotted}}

\newcommand{\Stwo}{{\small$p=2$}}
\newcommand{\Sthree}{{\small$p=3$}}

\newcommand{\ANS}{{\small ANS}}

\newcommand{\SStwo}{{\small SS}}
\newcommand{\SSANStwo}{{\small SS$_{\text{ANS}}$}}
\newcommand{\reference}{{\small Ref.}}

\begin{document}

\title{A simple and effective method based on strain projections to alleviate locking in isogeometric solid shells}

\author{Pablo Antolin \and Josef Kiendl \and Marco Pingaro \and Alessandro Reali}

\institute{P. Antolin (\Letter) \at
							\'Ecole Polytechnique F\'ed\'erale de Lausanne\\
              Institute of Mathematics\\
		          Lausanne, Switzerland.\\
              \email{pablo.antolin@epfl.ch}
							\and
							J.\ Kiendl \at
              Norwegian University of Science and Technology\\
              Department of Marine Technology\\
							Trondheim, Norway.
							\and
							M.\ Pingaro \at
              Universit\`a degli Studi di Roma ``La Sapienza''\\
							Dipartimento di Ingegneria Strutturale e Geotecnica\\
		          Roma, Italy.
							\and
							A.\ Reali \at
              Universit\`a degli Studi di Pavia\\
              Department of Civil Engineering and Architecture\\
							Pavia, Italy. \at
              Istituto di Matematica Applicata e\\ Tecnologie Informatiche ``E. Magenes'' (CNR)\\
							Pavia, Italy.
							}

\date{October 7th, 2019}

\maketitle

\begin{abstract}
In this work, we focus on the family of shell formulations referred to as ``solid shells'', where the simulation of shell-type structures is performed by means of a mesh of 3D solid elements, with typically only one element through the thickness. 
We propose a novel approach for alleviating the various locking phenomena, which typically appear in thin structures, based on the projection of strains onto discontinuous coarser polynomial spaces defined at element level. In particular, we present and investigate two different formulations based on this approach. Several numerical experiments prove the very good performance of both formulations. The main advantages of the presented approach compared to existing solid shell formulations are its simplicity and numerical efficiency. 
\end{abstract}

\section{Introduction \label{sec:introduction}}
Shell structures are ubiquitous in various engineering disciplines and there exists a huge variety of shell elements for analyzing them within finite element methods. As a major classification of shell elements, one can distinguish between (bivariate) shells and solid shells. While the former are formulated on bivariate domains and derived from shell theories, solid shells are trivariate solid elements and their formulation is based on 3D continuum theory.
Whereas standard solid elements applied to thin structures typically require excessively fine meshes with several elements through the  thickness to avoid geometric locking phenomena, in particular shear, membrane and curvature-thickness locking, solid shell elements are designed such that accurate analysis can be obtained with only one element in the thickness direction. Essentially, solid shells are solid elements enhanced by certain anti-locking techniques, like, e.g., $\bar{\boldsymbol{B}}-$ and  $\bar{\boldsymbol{F}}-$formulations, assumed natural strains (ANS), or enhanced assumed strains (EAS). 
The advantages of solid shells compared to bivariate shell elements are, among others:
(\textbullet) generally simpler formulation and implementation (standard solid theory vs.\ shell theories); 
(\textbullet) avoidance of rotational degrees of freedom, which are necessary in most classical shell elements; 
(\textbullet) straight-forward use of nonlinear constitutive models, which are generally derived in the context of 3D solids; 
(\textbullet) higher accuracy when three-dimensional stress states are important locally, e.g., for double-sided contact in sheet metal forming simulations.
Solid shells are well established in classical finite element analysis \cite{Roehl_1996,Hauptmann_1998,Sze_2000,Alves_2003,Valente_2004,Cardoso2008,Korelc_2010,Schwarze2011}.

In this paper we make use of the peculiar features of Isogeometric Analysis (IGA) (see, e.g., the monograph \cite{Cottrell2009} or the recent special issue \cite{CMAME-IGA2017}), which have been shown to have a great potential in particular for structural analysis.
So far, most of the IGA structural formulations proposed in the literature have been developed in the framework of bivariate shells (see, e.g., \cite{kiendl_isogeometric_2009,benson_large_2011,thai_static_2012,echter_hierarchic_2013,benson_blended_2013,hosseini_isogeometric_2013,kiendl_isogeometric_2015,oesterle_shear_2016,oesterle_hierarchic_2017,ambati_isogeometric_2018,leonetti_simplified_2019} and references therein), while only a few papers deal with IGA solid shell elements \cite{bouclier_efficient_2013,hosseini_isogeometric_2013,Hosseini2014,Caseiro2014,Caseiro2015,leonetti_efficient_2018}.

In this context, it is well-known that in IGA, locking is generally less pronounced due to the higher-order nature of the underlying function spaces, and it can be easily reduced to a practically insignificant level by simple order elevation. However, it is also well-known that higher-order IGA with standard quadrature rules can become numerically very costly and, therefore, there is a high interest in developing isogeometric solid shell elements which combine high-order accuracy with the efficiency of low-order approximations, i.e., quadratic elements. 
In \cite{Caseiro2014}, the concept of assumed natural strains (ANS) \cite{hughes_finite_1981,dvorkin_continuum_1984} was firstly applied to isogeometric solid shells. 
The general idea of the ANS method is to replace the strain components which cause locking by an ``assumed'' strain field. This approach can be summarized as follows: The compatible strains are evaluated at the so-called tying points instead of the integration points, where the tying points correspond to points of a reduced integration rule; an assumed strain field is then extrapolated from these tying points at element level; the assumed strain field is finally used in the weak form and integrated at the standard integration points.
Different sets of tying points are used for different strain components, which makes the implementation a bit cumbersome and also this affects numerical efficiency, due to the increased number of shape function evaluations (at the different sets of tying points and at the standard integration points), and, consequently, the memory requirements.
Another important aspect is that this method requires the general element formulation to be set in a curvilinear or local Cartesian coordinate frame for being able to separate the strains into in-plane and out-of-plane components. \par
In the present paper we propose two novel solid shell formulations, where locking is counteracted by projecting those strain components which cause locking onto coarser polynomial spaces at element level through local $L^2$ projections.
The first formulation is inspired by the ANS method in the sense that the projection spaces correspond to the spaces of the assumed strain fields in the ANS method, with the consequence that different projection operators are used for different strain components.
Its advantage, compared to the ANS method is that no tying points are necessary and shape functions are evaluated only at the standard integration points, which enhances the efficiency of the method.
Secondly, we explore a simplified formulation where the same projection is applied to all strain components.
In this case, the whole formulation and implementation becomes much more efficient.
In fact, this formulation does not even need a local coordinate system (as it is the case for many solid shell formulations), which makes its implementation into existing standard solid elements straight forward and particularly easy.
Numerical studies on several benchmark examples show that both formulations perform well, showing the same level of accuracy as the ANS formulation \cite{Caseiro2014}.
In this paper, we use NURBS basis functions for the discretization, however, the proposed formulations can be equally applied to standard Lagrangian finite elements. \par

\section{Standard solid and solid shell formulations}
The formulations presented in this work are based on classical 3D linear elasticity.
Thus, assuming, for the sake of simplicity, a combination of Dirichlet and homogeneous Neumann boundary conditions, the problem's weak form can be written as
\begin{align} \label{eq:elasticity_problem}
  \int_\Omega \delta\bm\varepsilon:\bm\sigma \,\diff\Omega =   \int_\Omega \delta\bm u\cdot\bm f \,\diff\Omega \, ,
\end{align}
where $\Omega$ denotes the domain occupied by the elastic body; $\bm u$ is the elastic displacement vector field (assumed to satisfy the prescribed Dirichlet boundary conditions); $\bm\sigma = \mathbb C : \bm\varepsilon$ is the stress tensor, being $\mathbb C$ the elasticity tensor, and $\bm\varepsilon = \nabla^S \bm u$ is the strain tensor ($\nabla^S$ denotes the symmetric gradient operator); and $\bm f$ is the external load vector field.
Finally, $\delta \bm u$ is a virtual displacement vector field (that satisfies homogenous Dirichlet boundary conditions where displacements are prescribed) and $\delta \bm\varepsilon$ is the virtual strain tensor field.

The elastic displacement $\bm u$ is approximated as
\begin{align}
  \bm u \approx \sum_{k=1}^{n} N_k(\xi,\eta,\zeta) \, \bm{u}_k \,,
  \label{u}
\end{align}
where $N_k\in\mathbb{R}$ are the NURBS basis functions \cite{Cottrell2009}, $n$ is the total number of functions, $(\xi,\eta,\zeta)$ are the coordinates in the parametric domain, and $\bm{u}_k\in\mathbb{R}^3$ are the control point displacements (i.e., the problem unknowns). Following a standard approach (see, e.g., \cite{Hughes1987}), we obtain the classical formulation for the element stiffness matrix to be
\begin{align}
  \bm k_e &= \int_{\Omega_e} \bm B^\top \bm D \bm B \, \diff\Omega_e\,,
  \label{k_e}
\end{align}
where $\Omega_e$ is the element domain, $\bm B$ is the strain-displacement matrix, and $\bm D$ is the material matrix. 
The strain-displacement matrix can generally be constructed as
\begin{align}
 \bm B = 
 \left(\begin{array}{cccc} \bm B_1 & \bm B_2  & \ldots  & \bm B_{n_e} \end{array} \right)\,,
 \label{B}
\end{align}
being $n_e$ the number of shape functions per element, with the submatrices $\bm B_k$ defined as
\begin{align}
 \bm B_k = 
 \left(\begin{array}{ccc} 
    N_{k,x} & 0 & 0 \\ 0 & N_{k,y} & 0 \\ 0 & 0 & N_{k,z} \\
    N_{k,y} & N_{k,x} & 0 \\ N_{k,z} & 0 & N_{k,x} \\ 0 & N_{k,z} & N_{k,y}
     \end{array} \right)\,,
     \label{Bca}
\end{align}
where the comma subscript indicates a partial derivative, e.g., $N_{k,x} = \partial N_k/\partial x$.
Thus, the strain tensor can be approximated as
\begin{align}
 \underline{\bm\eps} \approx \sum_{k=1}^{n} \bm B_k \, \bm{u}_k \,,
 \label{strain_cart}
\end{align}
where $\underline{\bm\eps}$ is the Voigt representation of the strain tensor in 
the Cartesian coordinate system, i.e.,
$\underline{\bm\eps} = [\eps_{xx},\allowbreak\,\eps_{yy},\allowbreak\,\eps_{zz},\allowbreak\,2\eps_{xy},\allowbreak2\,\eps_{xz},\allowbreak\,2\eps_{yz}]^\top$.
\par

\begin{figure}[!ht]
\begin{center}
  \tikzsetnextfilename{shell_kinematics}\input{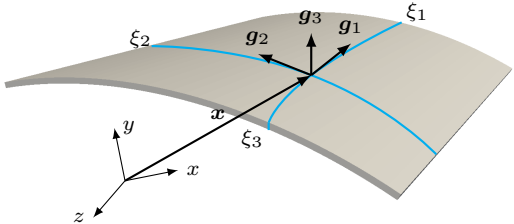}
  \caption{Solid shell curvilinear coordinates and its associated covariant basis.}\label{fig:shell_kinematics}
\end{center}
\end{figure}

In many solid shell formulations it is necessary to express the strains and, accordingly, the strain-dis\-place\-ment matrix in a curvilinear coordinate system aligned with the shell's geometry in order to separate the strains into in-plane and out-of-plane components. Typically, the geometry is modeled such that the first two coordinates $(\xi,\eta)$ correspond to the in-plane directions of the shell and the third one $(\zeta)$ to the thickness direction. 
We can then compute the curvilinear base vectors $\bm g_i$ as
\begin{align}
  \bm g_i &= \frac{\partial\bm x}{\partial\xi_i} = \sum_{k=1}^{n} \frac{\partial N_k}{\partial \xi_i} \bm x_k\,,
  \quad i=\lbrace1,\,2,\,3\rbrace\,,
  \label{gi}
\end{align}
where we used $(\xi_1,\xi_2,\xi_3)=(\xi,\eta,\zeta)$ for a shorter notation (see Figure \ref{fig:shell_kinematics}), and $\bm{x}_k$ are the geometry control point coordinates.

Using the basis \eqref{gi}, we can compute the strain-displacement matrix, referred to the curvilinear system and in a row-wise way, as
\begin{align}
 \widetilde{\bm B}_k = 
 \left(\begin{array}{c}   
     N_{k,\xi}\,\bm g_1^\top  \\  N_{k,\eta}\,\bm g_2^\top  \\  N_{k,\zeta}\,\bm g_3^\top  \\
     N_{k,\xi}\,\bm g_2^\top + N_{k,\eta}\,\bm g_1^\top  \\  
     N_{k,\xi}\,\bm g_3^\top + N_{k,\zeta}\,\bm g_1^\top  \\  
     N_{k,\eta}\,\bm g_3^\top + N_{k,\zeta}\,\bm g_2^\top  
 \end{array} \right)  \,.
 \label{Bcu}
\end{align}
Analogously as for the Cartesian system \eqref{strain_cart}, the 
curvilinear (covariant) strain components $\underline{\widetilde{\bm\eps}} = [\tilde\eps^{xx},\allowbreak\,\tilde\eps^{yy},\allowbreak\,\tilde\eps^{zz},\allowbreak\,2\tilde\eps^{xy},\allowbreak\,2\tilde\eps^{xz},\allowbreak\,2\tilde\eps^{yz}]^\top$ can be expressed as
\begin{align}
  \underline{\widetilde{\bm\eps}} \approx \sum_{k=1}^{n} \widetilde{\bm B}_k \, \bm{u}_k \,,
 \label{strain_cov}
\end{align}
where, in the same way as in \eqref{B}, the submatrices $\widetilde{\bm B}_k$ can be gathered as 
\begin{align} \label{eq:full_B}
 \widetilde{\bm B} = 
 \left(\begin{array}{cccc} \widetilde{\bm B}_1 & \widetilde{\bm B}_2  & \ldots  & \widetilde{\bm B}_{n_e} \end{array} \right).
\end{align}
It should be noted that in \eqref{Bcu} the derivatives with respect to the natural NURBS coordinates $(\xi,\eta,\zeta)$ are used.
Also note that, in contrast to the classical matrix $\bm B$, the curvilinear matrix $\widetilde{\bm B}$ is, in general, fully populated, as it can be observed comparing Equations \eqref{Bca} and  \eqref{Bcu}.
 
For consistency, also the Cartesian material matrix $\bm D$ needs to be expressed in the curvilinear system. This is obtained via the transformation matrix $\bm R$ (see, for instance, \cite{bouclier_efficient_2013}):
\begin{equation}\label{Rmat}
\begin{aligned}
 \bm R &=
 \left[\begin{matrix}
    J_{11}^2 & J_{21}^2 & J_{31}^2 & J_{11}J_{21}\\
    J_{12}^2 & J_{22}^2 & J_{32}^2 & J_{12}J_{22}\\
    J_{13}^2 & J_{23}^2 & J_{33}^2 & J_{13}J_{23}\\
    2J_{11}J_{12} & 2J_{22}J_{21} & 2J_{31}J_{32} & J_{11}J_{22}+J_{21}J_{12}\\
    2J_{11}J_{13} & 2J_{21}J_{23} & 2J_{31}J_{33} & J_{11}J_{23}+J_{21}J_{13}\\
    2J_{12}J_{13} & 2J_{22}J_{23} & 2J_{32}J_{33} & J_{12}J_{23}+J_{22}J_{13}
 \end{matrix}\right.\\
 &\qquad\qquad\quad
 \left.\begin{matrix}
 J_{11}J_{31} & J_{21}J_{31} \\
 J_{12}J_{32} & J_{22}J_{32} \\
 J_{13}J_{33} & J_{23}J_{33} \\
 J_{11}J_{32}+J_{31}J_{12} & J_{21}J_{32}+J_{31}J_{22} \\
 J_{11}J_{33}+J_{31}J_{13} & J_{21}J_{33}+J_{31}J_{23} \\
 J_{12}J_{33}+J_{32}J_{13} & J_{22}J_{33}+J_{32}J_{23}
 \end{matrix}\right]
\end{aligned}
\end{equation}
where $J_{ij}$ are the components of the Jacobian matrix:
\begin{align}
  \bm J =
  \left(\begin{array}{ccc} \frac{\partial x}{\partial\xi} & \frac{\partial x}{\partial\eta} & \frac{\partial x}{\partial\zeta} \\
   \frac{\partial y}{\partial\xi} & \frac{\partial y}{\partial\eta} & \frac{\partial y}{\partial\zeta}  \\ 
   \frac{\partial z}{\partial\xi} & \frac{\partial z}{\partial\eta} & \frac{\partial z}{\partial\zeta}   \end{array} \right) 
  =  \left(\begin{array}{ccc} \bm g_1 & \bm g_2 & \bm g_3 \end{array} \right) .
\end{align}
The material matrix in the curvilinear system is then obtained as
\begin{align}
  \widetilde{\bm D}  =  \bm R^{-T}\bm D \bm R^{-1}\,,
\end{align}
and the element stiffness matrix can be computed as
\begin{align}
  \bm k_e &= \int_{\Omega_e} \widetilde{\bm B}^\top \widetilde{\bm D} \widetilde{\bm B} \,\diff\Omega_e.
  \label{k_e_curv}
\end{align}
Note that, in the same way, $\widetilde{\bm B}$ can be alternatively calculated as $\widetilde{\bm B} = \bm R\bm B$.

\section{Solid shell locking alleviation} \label{sec:locking}

In this section, we present two different approaches to alleviate locking effects through local modifications of the strain-displacement matrices, for both the curvilinear \eqref{k_e_curv} and the Cartesian \eqref{k_e} formulations.

To this end, we propose a procedure based on the projection of the strain components causing locking onto coarser polynomial spaces.
In fact, in the same way as in other problems, like, e.g.,\ quasi-incompressible elasticity, locking phenomena arise as a consequence of the excess of constraints on the numerical solution.
In order to reduce these constraints, the affected part of the elastic energy is interpolated by means of coarser polynomial spaces.
In particular, for the case of solid shells, we interpolate the different components of the strain tensor using element-wise lower order polynomial spaces that are discontinuous across elements.

\subsection{ANS-inspired locking alleviation using local projections} \label{sec:ANS-inspired}

Let us first introduce the strain tensor expressed in the covariant basis:
\begin{align}
  \bm{\eps} = \sum^{3}_{i,j=1} \tilde\eps^{ij} \bm{g}_i\otimes\bm{g}_j\,.
\end{align}
In order to alleviate the possible locking phenomena, the strain tensor $\bm{\eps}$ is replaced with its modified version $\overline{\bm{\eps}}$ that results from a local $L^2$ projection of each covariant component of $\bm{\eps}$ onto different coarser polynomial spaces defined for each parametric element $\hat\Omega_e$.
Accordingly, each component $\tilde\eps^{ij}$ of $\bm{\eps}$ is substituted by a new component $\overline{\eps}^{ij}$ of $\overline{\bm{\eps}}$, obtained through its $L^2$ projection onto a reduced order space.
Thus, the modified strain tensor reads
\begin{align}
  \overline{\bm{\eps}} = \sum^{3}_{i,j=1} \overline{\eps}^{ij} \bm{g}_i\otimes\bm{g}_j\,,\quad\text{with }\overline{\eps}^{ij} = \Pi^{\left(i,j\right)}(\tilde\eps^{ij})\,,
\end{align}
where $\Pi^{\left(i,j\right)}$ is the element-wise $L^2$ projection operator onto the space $\mathbb{Q}^{(i,j)}(\hat{\Omega}_e)$ for each $ij$ strain component.

Thus, in each single parametric element $\hat{\Omega}_e$, the $L^2$ projection can be implicitly expressed as:
\begin{align}\label{eq:projectors2}
\begin{split}
  \int_{\hat{\Omega}_e} \tilde\eps^{ij} \theta_h\dOp_e &= \int_{\hat{\Omega}_e} \overline{\eps}^{ij}\theta_h\dOp_e\,, \\
  &\forall\theta_h, \overline{\eps}^{ij} \in\mathbb{Q}^{(i,j)}(\hat{\Omega}_e),\,\, \tilde\eps^{ij} \in L^2(\hat{\Omega}_e)\,.
\end{split}
\end{align}

Assuming that the problem solution is discretized with the same degree $p$ along the three parametric directions, and inspired by the ANS method proposed for isogeometric analysis in \cite{Caseiro2014,Caseiro2015}, the different strain components are then treated as follows:
\begin{itemize}
  \item $\tilde\eps^{11}$ and $\tilde\eps^{13}$ are projected onto the local element space $\mathbb{Q}^{(1,1)}=\mathbb{Q}_{p-1,p,p}(\hat\Omega_e)$;
  \item $\tilde\eps^{22}$ and $\tilde\eps^{23}$ are projected onto the local element space $\mathbb{Q}^{(2,2)}=\mathbb{Q}_{p,p-1,p}(\hat\Omega_e)$;
  \item $\tilde\eps^{12}$ is projected onto the local element space $\mathbb{Q}^{(1,2)}=\mathbb{Q}_{p-1,p-1,p}(\hat\Omega_e)$;
  \item $\tilde\eps^{33}$ remains unprojected, i.e., $\overline{\eps}_{33} = \tilde\eps_{33}$.
\end{itemize}
$\mathbb{Q}_{q,r,s}(\hat\Omega_e)$ is the space of polynomials of degrees $\leq(q,r,s)$, along the three parametric directions of the parametric domain element $\hat\Omega_e$.

The above-introduced polynomial spaces are local to each element and discontinuous across different elements.
Therefore, the projections of the different strain components can be computed at each element independently from the others, making this operation computationally inexpensive and embarrassingly parallel.

Thus, the element-wise projection can be explicitly written in matrix form as:
\begin{align}
  \overline{\Matrix{f}}^{ij} = \mathbb{P}^{(i,j)}\, \tilde{\Matrix{f}}^{ij}\,,
\end{align}
where the column vectors $\tilde{\Matrix{f}}^{ij}\in\mathbb{R}^{n_q}$ and $\overline{\Matrix{f}}^{ij}\in\mathbb{R}^{n_q}$ are the values of $\tilde\eps^{ij}$ and $\overline{\eps}^{ij}$, respectively, evaluated at the $n_q$ quadrature points of a single element $\hat\Omega_e$, while $\mathbb{P}^{(i,j)}\in\mathbb{R}^{n_q\times n_q}$ is the linear projection operator expressed in matrix form.
Due to the fact that the projections are performed in the parametric domain, and the same projection spaces are chosen for all elements, the operators $\mathbb{P}^{(i,j)}$ are constant from element to element.
In Appendix \ref{sec:proj_coeffs} we provide closed-form expressions of the matrices $\mathbb{P}^{(i,j)}$, for degrees $p=1$ and $p=2$, ready to be used in isogeometric or finite element analysis codes.

By means of the above defined projections, the modified strain can be represented as
\begin{align}
  \underline{\overline{\bm\eps}} \approx \sum_{k=1}^{n} \overline{\bm B}_k \, \bm{u}_k \,,
 \label{proj_strain_cov}
\end{align}
where $\overline{\bm B}_k$ is written row-wise as
\begin{align}
 \overline{\bm B}_k = 
 \left(\begin{array}{c}   
     \Pi^{(1,1)}\left(N_{k,\xi}\,\bm g_1^\top\right)  \\
     \Pi^{(2,2)}\left(N_{k,\eta}\,\bm g_2^\top\right)  \\
     N_{k,\zeta}\,\bm g_3^\top  \\
     \Pi^{(1,2)}\left(N_{k,\xi}\,\bm g_2^\top + N_{k,\eta}\,\bm g_1^\top\right)  \\  
     \Pi^{(1,1)}\left(N_{k,\xi}\,\bm g_3^\top + N_{k,\zeta}\,\bm g_1^\top\right)  \\  
     \Pi^{(2,2)}\left(N_{k,\eta}\,\bm g_3^\top + N_{k,\zeta}\,\bm g_2^\top\right)  
 \end{array} \right)  \,.
 \label{Bcu_proj}
\end{align}
In the same way as for  \eqref{eq:full_B}, we can define 
\begin{align}
 \overline{\bm B} = 
 \left(\begin{array}{cccc} \overline{\bm B}_1 & \overline{\bm B}_2  & \ldots  & \overline{\bm B}_{n_e} \end{array} \right),
\end{align}
such that the element  stiffness matrix is computed as
\begin{align}
  \overline{\bm k}_e &= \int_{\Omega_e} \overline{\bm B}^\top \widetilde{\bm D} \overline{\bm B} \,\diff\Omega_e.
  \label{k_e_curv_proj}
\end{align}

\subsection{Simplified Cartesian locking alleviation using local projections} \label{sec:cartesian}

In this work we also explore the possibility of projecting all strain components, including $\tilde\varepsilon^{33}$,
onto the reduced polynomial space $\mathbb{Q}_{p-1,p-1,p}(\hat{\Omega}_e)$.
In such a case, all strain components are projected directly using their Cartesian version,
and their transformation to curvilinear components is no longer needed, which renders the operation much simpler.

Thus, the new strain-displacement matrix $\widehat{\bm B}_k$ is computed as:
\begin{align}
 \widehat{\bm B}_k = 
 \left(\begin{array}{ccc} 
   \widehat{N}_{k,x} & 0 & 0 \\ 0 & \widehat{N}_{k,y} & 0 \\ 0 & 0 & \widehat{N}_{k,z} \\
   \widehat{N}_{k,y} & \widehat{N}_{k,x} & 0 \\ \widehat{N}_{k,z} & 0 & \widehat{N}_{k,x} \\ 0 & \widehat{N}_{k,z} & \widehat{N}_{k,y}
     \end{array} \right)\,,
\end{align}
where $\widehat{N}_{k,x_i}$ are simply the projected basis function derivatives, i.e.,
\begin{align} \label{B_3}
  \widehat{N}&_{k,x_i} = \Pi^{(1,2)}(N_{k,x_i})\,,
\end{align}
for $i=\lbrace1,2,3\rbrace$ and $(x_1,\,x_2,\,x_3)=(x,\,y,\,z)$.
Then, as before, the local stiffness matrix is simply computed as:
\begin{align}
  \widehat{\bm k}_e &= \int_{\Omega_e} \widehat{\bm B}^\top \bm D \widehat{\bm B} \, \diff\Omega_e\,.
  \label{k_e_3}
\end{align}

This methodology presents a fundamental advantage with respect to the projections presented in the previous section: It requires the use of Cartesian coordinates only (notice the use of Cartesian coordinates of the involved derivatives in \eqref{B_3} and the use of the Cartesian version of tensor $\bm D$ in \eqref{k_e_3}), making it simpler and faster.
Moreover, it is worth mentioning that this approach differs from reduced integration techniques, which in some cases may lead to unphysical instabilities, while, on the basis of our numerical tests, the proposed method seems to be always stable.
As it will be shown in the numerical experiments gathered in Section \ref{sec:numerical}, the performance of this simplified technique is as good as the ANS-inspired one, presented in Section \ref{sec:ANS-inspired}, and even superior in some cases.

\section{Numerical examples} \label{sec:numerical}
In this section we present a series of numerical experiments, classically used to test shell and solid-shell formulation capabilities for alleviating locking effects, with the aim of illustrating the performance of the proposed solid-shell elements.

In particular we first present two classical beam tests, and  we then analyze the celebrated set of three benchmarks known as the ``shell obstacle course'', proposed by Belytschko et al.\ in \cite{Belytschko1985}.

In these five test cases we compare the performance of different quadratic solid shell elements, namely: The ANS-inspired version of our quadratic solid-shell element, described in Section \ref{sec:ANS-inspired} and labeled hereinafter as ``\ref{plt:line:solid_shell_ANS} \SSANStwo'';
the simplified quadratic solid-shell element, formulated in Cartesian coordinates and detailed in Section \ref{sec:cartesian}, denoted as ``\ref{plt:line:solid_shell} \SStwo'';
and, finally, the quadratic ANS element proposed in \cite{Caseiro2014} and labeled as ``\ref{plt:line:ANS} \ANS''.
Quadratic formulations are also compared with the standard cubic isogeometric element
(simply denoted as ``\ref{plt:line:s3} \Sthree''), known to show a good behavior
even in the presence of shear and membrane locking conditions.
Additionally, we also include in all the test cases the results corresponding to the standard isogeometric quadratic element, denoted as ``\ref{plt:line:s2} \Stwo''.

\subsection{Straight cantilever beam}
\begin{figure}
\begin{center}
  \tikzsetnextfilename{cantilever_geometry}\tdplotsetmaincoords{75}{45}
\begin{tikzpicture}[scale=0.7,tdplot_main_coords]
	
	\coordinate (O) at (0,0,0);
	\pgfmathsetmacro{\Px}{10}
	\pgfmathsetmacro{\Py}{2}
	\pgfmathsetmacro{\Pz}{0.5}
	
	\fill[fill=black!70!white,fill opacity=0.4] (0,-\Py,-\Pz) -- (0,2*\Py,-\Pz) -- (0,2*\Py,2*\Pz) -- (0,-\Py,2*\Pz) -- cycle;	
	
	\fill[fill=blue!50!white,fill opacity=0.3] (0,0,0) -- (\Px,0,0) -- (\Px,\Py,0) -- (0,\Py,0) -- cycle;
	\fill[fill=blue!50!white,fill opacity=0.3] (0,0,\Pz) -- (\Px,0,\Pz) -- (\Px,\Py,\Pz) -- (0,\Py,\Pz) -- cycle;
	\fill[fill=blue!50!white,fill opacity=0.3] (\Px,0,0) -- (\Px,\Py,0) -- (\Px,\Py,\Pz) -- (\Px,0,\Pz) -- cycle;	
	\fill[fill=blue!50!white,fill opacity=0.3] (0,0,0) -- (0,\Py,0) -- (0,\Py,\Pz) -- (0,0,\Pz) -- cycle;
	\fill[fill=blue!50!white,fill opacity=0.3] (0,0,0) -- (\Px,0,0) -- (\Px,0,\Pz) -- (0,0,\Pz) -- cycle;
	\fill[fill=blue!50!white,fill opacity=0.3] (0,\Py,0) -- (\Px,\Py,0) -- (\Px,\Py,\Pz) -- (0,\Py,\Pz) -- cycle;
	
	\draw[line width=0.7pt, draw=black] (0,0,0) -- (\Px,0,0) -- (\Px,\Py,0); 
	\draw[line width=0.7pt, draw=black, dashed] (\Px,\Py,0) -- (0,\Py,0) -- (0,0,0);
	\draw[line width=0.7pt, draw=black] (0,0,\Pz) -- (\Px,0,\Pz) -- (\Px,\Py,\Pz) -- (0,\Py,\Pz) -- cycle;
	\draw[line width=0.7pt, draw=black] (0,0,0) -- (0,0,\Pz);
	\draw[line width=0.7pt, draw=black, dashed] (0,\Py,0) -- (0,\Py,\Pz);
	\draw[line width=0.7pt, draw=black] (\Px,0,0) -- (\Px,0,\Pz);
	\draw[line width=0.7pt, draw=black] (\Px,\Py,0) -- (\Px,\Py,\Pz);
	\draw[line width=0.8pt,->] (0,-\Py,-2) -- (1,-\Py,-2) node[anchor=north east]{$x$};
	\draw[line width=0.8pt,->] (0,-\Py,-2) -- (0,-\Py+1,-2) node[anchor=south west]{$y$};
	\draw[line width=0.8pt,->] (0,-\Py,-2) -- (0,-\Py,-1) node[anchor=north east]{$z$};	
	
	\begin{scope}[->, ultra thick]
    \draw[->,line width=1.3pt,red] (\Px,0,\Pz+1.5)--(\Px,0,\Pz);
    \draw[->,line width=1.3pt,red] (\Px,\Py/3,\Pz+1.5)--(\Px,\Py/3,\Pz);
    \draw[->,line width=1.3pt,red] (\Px,2*\Py/3,\Pz+1.5)--(\Px,2*\Py/3,\Pz);
    \draw[->,line width=1.3pt,red] (\Px,\Py,\Pz+1.5)--(\Px,\Py,\Pz);
    \end{scope}
    \draw [line width=1.3pt,color=red] (\Px,0,\Pz+1.5) -- (\Px,\Py,\Pz+1.5);

	\draw [line width=0.5,color=black, dashed] (O) -- (0,0, -\Pz*2); 		
	\draw [line width=0.5,color=black, dashed] (\Px,0,0) -- (\Px, 0, -\Pz*2);	
	\draw [<->,color=black] (0,0, -\Pz*2) -- (\Px, 0, -\Pz*2) node[black,midway,above] {$L$};
	
	\draw [line width=0.5,color=black, dashed] (\Px,0,0) -- (\Px+2,0, 0);
	\draw [line width=0.5,color=black, dashed] (\Px,\Py,0) -- (\Px+2, \Py, 0);
	\draw [<->,color=black] (\Px+2,0, 0) -- (\Px+2, \Py, 0) node[black,midway,below] {$w$};
	
	\draw [line width=0.5,color=black, dashed] (\Px,\Py,0) -- (\Px,\Py+2,0);
	\draw [line width=0.5,color=black, dashed] (\Px,\Py,\Pz) -- (\Px, \Py+2, \Pz);
	\draw [<->,color=black] (\Px,\Py+2,0) -- (\Px, \Py+2, \Pz) node[black,sloped,midway,right, rotate=-90] {$t$};
	
	\draw [color=red] (\Px,\Py/2,\Pz+2.2) node[anchor=north, rotate=0] {$F$};
\end{tikzpicture}
  \caption{Straight cantilever beam: Problem description.
	One face is clamped whereas a distributed load (with resultant $F$) is applied along the top edge
	of the opposite face. Beam deflection is measured at the bottom edge of the free end.
  $E=1000$, $\nu=0$, $L=100$, and $w=1$. Different slendernesses $L/t$ are considered.}\label{fig:cantilever_straight_geom}
\end{center}
\end{figure}
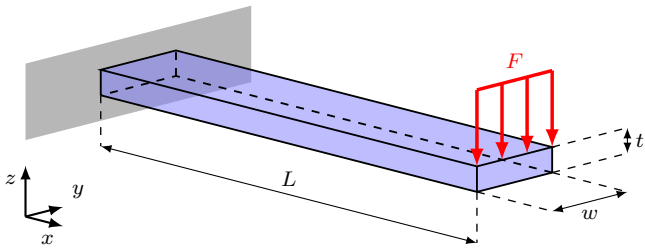

This first example is a straight cantilever beam, clamped at one face and subjected to a distributed load, 
with resultant $F$, along the top edge of the opposite face (see Figure \ref{fig:cantilever_straight_geom}).
A Young modulus $E=1000$ and a Poisson ratio $\nu=0$ are assumed, while the beam length and width are $L=100$ and $w=1$, respectively. The beam thickness is indicated by $t$.
In all  numerical tests a single element is used for the beam cross section, while a variable number of elements along the longitudinal direction are considered.

Even if the nature of the formulation proposed in this paper is three dimensional, under the geometrical setting, boundary, loading, and material conditions
described above, the model behaves as a 2D beam.
This test helps in evaluating the effect of shear locking, isolated from other possible effects, for
all the considered discretizations.

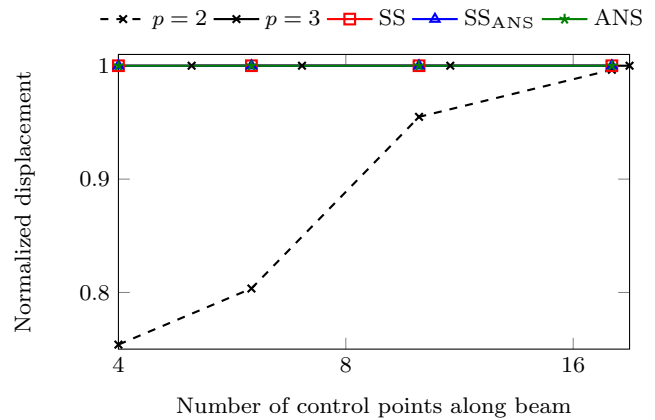
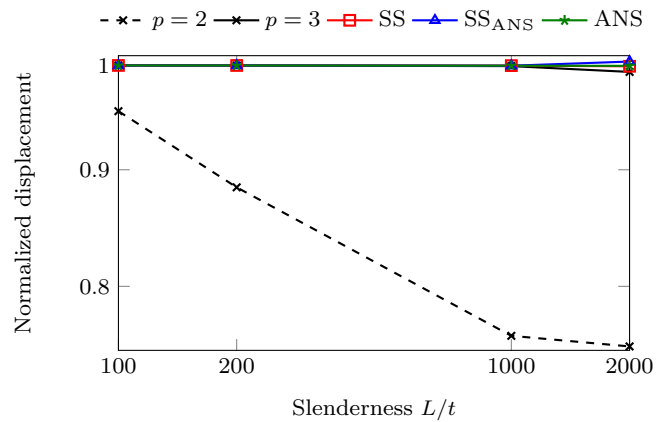
\begin{figure}
\centering
\subfigure[Beam vertical deflection as a function of the number of control points along beam length  (slenderness $L/t=100$).\label{fig:cantilever_convergence_lin}]{\tikzsetnextfilename{straight_beam_convergence}\begin{tikzpicture}
 \begin{semilogxaxis}[basic_paper_plots_style,
  xlabel={Number of control points along beam},
  ylabel={Normalized displacement},
  ymin=0.75, ymax=1.01,
  xmin=4, xmax=19,
  xtick={4, 8, 16},
  xticklabels={$4$, $8$, $16$},
  ytick={0.8, 0.9, 1},
  yticklabels={$0.8$, $0.9$, $1$},
  scale=0.45,
  legend columns=5
 ]
 %
 %
 \addplot[solid 2]
 table [x=num_ctr_pts_side_2, y=S2]
 {tikz_figures/straight_beam_Rt_100.res};\label{plt:line:s2}
 \addlegendentry{\Stwo}
 \addplot[solid 3]
 table [x=num_ctr_pts_side_3, y=S3]
 {tikz_figures/straight_beam_Rt_100.res};\label{plt:line:s3}
 \addlegendentry{\Sthree}
 \addplot[ss]
 table [x=num_ctr_pts_side_2, y=SS2]
 {tikz_figures/straight_beam_Rt_100.res};\label{plt:line:solid_shell}
 \addlegendentry{\SStwo}
 \addplot[ss ANS-like]
 table [x=num_ctr_pts_side_2, y=SS2-ANS-LIKE]
 {tikz_figures/straight_beam_Rt_100.res};\label{plt:line:solid_shell_ANS}
 \addlegendentry{\SSANStwo}
 \addplot[ANS]
 table [x=num_ctr_pts_side_2, y=ANS2]
 {tikz_figures/straight_beam_Rt_100.res};\label{plt:line:ANS}
 \addlegendentry{\ANS}
 \end{semilogxaxis}
\end{tikzpicture}}
\subfigure[Beam vertical deflection as a function of slenderness $L/t$ (8 elements along the beam's length).\label{fig:cantilever_slen_lin}]{\tikzsetnextfilename{straight_beam_slenderness}\begin{tikzpicture}
 \begin{loglogaxis}[basic_paper_plots_style,
  xlabel={Slenderness $L/t$},
  ylabel={Normalized displacement},
  ymin=0.75, ymax=1.01,
  xmin=100, xmax=2000,
  ytick={0.8, 0.9, 1},
  yticklabels={$0.8$, $0.9$, $1$},
  xtick={100, 200, 1000, 2000},
  xticklabels={$100$, $200$, $1000$, $2000$},
  scale=0.45,
  legend columns=5
 ]
 %
 %
 \addplot[solid 2]
 table [x=slenderness, y=S2]
 {tikz_figures/straight_beam_nel_8.res};
 \addlegendentry{\Stwo}
 \addplot[solid 3]
 table [x=slenderness, y=S3]
 {tikz_figures/straight_beam_nel_8.res};
 \addlegendentry{\Sthree}
 \addplot[ss]
 table [x=slenderness, y=SS2]
 {tikz_figures/straight_beam_nel_8.res};
 \addlegendentry{\SStwo}
 \addplot[ss ANS-like]
 table [x=slenderness, y=SS2-ANS-LIKE]
 {tikz_figures/straight_beam_nel_8.res};
 \addlegendentry{\SSANStwo}
 \addplot[ANS]
 table [x=slenderness, y=ANS2]
 {tikz_figures/straight_beam_nel_8.res};
 \addlegendentry{\ANS}
 \end{loglogaxis}
\end{tikzpicture}}
\caption{Straight cantilever beam: Normalized vertical deflection at the tip for different discretizations and slendernesses.}
\end{figure}

Considering a high slenderness value $L/t=100$, the normalized deflection of the beam tip is reported in Figure \ref{fig:cantilever_convergence_lin}
for different numbers of control points along the beam length.
As it can be seen, all considered solid-shell elements are able to capture the exact
solution, even for the coarsest considered mesh (this is not the case for standard quadratic elements).

On the other hand, in Figure \ref{fig:cantilever_slen_lin} we report the results obtained in the case of a fixed discretization (8 elements along the beam's length) and different (high) values of the slenderness $L/t$.
As in the previous case, except for the standard quadratic element, all discretizations
are able to capture the correct results, even for quite severe slendernesses.

\begin{figure}
  \centering
  \subfigure[$30^\circ$ in-plane distortion.\label{fig:cantilever_distorted_lin_30}]{\tikzsetnextfilename{straight_beam_slenderness_distorted_30}\begin{tikzpicture}
 \begin{semilogxaxis}[basic_paper_plots_style,
  xlabel={Slenderness $L/t$},
  ylabel={Normalized displacement},
  ymin=0.4, ymax=1.01,
  ytick={0.4, 0.5, 0.6, 0.7, 0.8, 0.9, 1.0},
  yticklabels={$0.4$, $0.5$, $0.6$, $0.7$, $0.8$, $0.9$, $1$},
  xmin=100, xmax=2000,
  xtick={100, 200, 1000, 2000},
  xticklabels={$100$, $200$, $1000$, $2000$},
	scale=0.45,
	legend columns=5,
 ]
 \addplot[solid 2]
 table [x=slenderness, y=S2]
 {tikz_figures/straight_beam_distorted_ang_30_nel_8.res};
 \addlegendentry{\Stwo}
 \addplot[solid 3]
 table [x=slenderness, y=S3]
 {tikz_figures/straight_beam_distorted_ang_30_nel_8.res};
 \addlegendentry{\Sthree}
 \addplot[ss ANS-like]
 table [x=slenderness, y=SS2-ANS-LIKE]
 {tikz_figures/straight_beam_distorted_ang_30_nel_8.res};
 \addlegendentry{\SSANStwo}
 \addplot[ss]
 table [x=slenderness, y=SS2]
 {tikz_figures/straight_beam_distorted_ang_30_nel_8.res};
 \addlegendentry{\SStwo}
 \addplot[ANS]
 table [x=slenderness, y=ANS2]
 {tikz_figures/straight_beam_distorted_ang_30_nel_8.res};
 \addlegendentry{\ANS}
 \end{semilogxaxis}
 \node at (3, 1.25) {\begin{tikzpicture}[scale = 0.6]
	
	\coordinate (o) at (0,0);	
	\pgfmathsetmacro{\L}{6}
	\pgfmathsetmacro{\H}{1}
	\pgfmathsetmacro{\l}{3}
	\pgfmathsetmacro{\ll}{\l-1}
	\pgfmathsetmacro{\tan}{0.5773502691896256}

	\pgfmathsetmacro{\e}{0.13}
	\draw[line width=0.7pt, draw=black] (0,0) -- (\L, 0) -- (\L, \H) -- (0, \H) -- cycle;
	\foreach \x in {1,2,...,\l}
	{
	  \pgfmathsetmacro{\aux}{0.5*\H*\tan/\l*(\l+1-\x)}
	  \pgfmathsetmacro{\a}{(\x-1)*0.5*\L/\l}
		\draw[line width=0.7pt, draw=black] (0.5*\L-\a+\aux,  0) -- (0.5*\L-\a-\aux, \H);

    \ifthenelse{\x>1}
    {
		\draw[line width=0.7pt, draw=black] (0.5*\L+\a+\aux,  0) -- (0.5*\L+\a-\aux, \H);
    }{}
	}
	



	\coordinate (P0) at (\L/2,\H/2);
	\coordinate (P1) at (\L/2,\H);
	\coordinate (P4) at (\L/2,2*\H);
	\coordinate (P2) at ($(P0) + (-0.5*\H*\tan,\H/2)$);
	\coordinate (P3) at ($(P0) + (-\H*\tan,\H)$);

  \draw[dashed] ($(P0)!0.1!(P4)$) -- ($(P0)!1.2!(P4)$);
  \draw[dashed] ($(P2)!0.1!(P3)$) -- ($(P2)!2.3!(P3)$);

  \draw ($(P0) + (120:\H*1.25)$) arc (120:90:\H*1.25) node [midway, above] {$30^\circ$};
	
  \coordinate (oo) at (0,1.25*\H);	
  \draw[stealth'-stealth'] ($(oo) + (0, 0.75)$) node[right] {\footnotesize $y$} -- (oo) -- ($(oo) + (0.75, 0)$) node[above] {\footnotesize $x$};


\end{tikzpicture}};
\end{tikzpicture}}
  \subfigure[$60^\circ$ in-plane distortion.\label{fig:cantilever_distorted_lin_60}]{\tikzsetnextfilename{straight_beam_slenderness_distorted_60}\begin{tikzpicture}
 \begin{semilogxaxis}[basic_paper_plots_style,
  xlabel={Slenderness $L/t$},
  ylabel={Normalized displacement},
  ymin=0.4, ymax=1.01,
  ytick={0.4, 0.5, 0.6, 0.7, 0.8, 0.9, 1.0},
  yticklabels={$0.4$, $0.5$, $0.6$, $0.7$, $0.8$, $0.9$, $1$},
  xmin=100, xmax=2000,
  xtick={100, 200, 1000, 2000},
  xticklabels={$100$, $200$, $1000$, $2000$},
	scale=0.45,
	legend columns=5,
 ]
 \addplot[solid 2]
 table [x=slenderness, y=S2]
 {tikz_figures/straight_beam_distorted_ang_60_nel_8.res};
 \addlegendentry{\Stwo}
 \addplot[solid 3]
 table [x=slenderness, y=S3]
 {tikz_figures/straight_beam_distorted_ang_60_nel_8.res};
 \addlegendentry{\Sthree}
 \addplot[ss ANS-like]
 table [x=slenderness, y=SS2-ANS-LIKE]
 {tikz_figures/straight_beam_distorted_ang_60_nel_8.res};
 \addlegendentry{\SSANStwo}
 \addplot[ss]
 table [x=slenderness, y=SS2]
 {tikz_figures/straight_beam_distorted_ang_60_nel_8.res};
 \addlegendentry{\SStwo}
 \addplot[ANS]
 table [x=slenderness, y=ANS2]
 {tikz_figures/straight_beam_distorted_ang_60_nel_8.res};
 \addlegendentry{\ANS}
 \end{semilogxaxis}
 \node at (3, 1.25) {\begin{tikzpicture}[scale = 0.6]
	
	\coordinate (o) at (0,0);	
	\pgfmathsetmacro{\L}{6}
	\pgfmathsetmacro{\H}{1}
	\pgfmathsetmacro{\l}{3}
	\pgfmathsetmacro{\ll}{\l-1}
	\pgfmathsetmacro{\tan}{1.7320508075688767}

	\pgfmathsetmacro{\e}{0.13}
	\draw[line width=0.7pt, draw=black] (0,0) -- (\L, 0) -- (\L, \H) -- (0, \H) -- cycle;
	\foreach \x in {1,2,...,\l}
	{
	  \pgfmathsetmacro{\aux}{0.5*\H*\tan/\l*(\l+1-\x)}
	  \pgfmathsetmacro{\a}{(\x-1)*0.5*\L/\l}
		\draw[line width=0.7pt, draw=black] (0.5*\L-\a+\aux,  0) -- (0.5*\L-\a-\aux, \H);

    \ifthenelse{\x>1}
    {
		\draw[line width=0.7pt, draw=black] (0.5*\L+\a+\aux,  0) -- (0.5*\L+\a-\aux, \H);
    }{}
	}
	



	\coordinate (P0) at (\L/2,\H/2);
	\coordinate (P1) at (\L/2,\H);
	\coordinate (P4) at (\L/2,2*\H);
	\coordinate (P2) at ($(P0) + (-0.5*\H*\tan,\H/2)$);
	\coordinate (P3) at ($(P0) + (-\H*\tan,\H)$);

  \draw[dashed] ($(P0)!0.1!(P4)$) -- ($(P0)!1.2!(P4)$);
  \draw[dashed] ($(P2)!0.1!(P3)$) -- ($(P2)!1.2!(P3)$);

  \draw ($(P0) + (150:\H*1.25)$) arc (150:90:\H*1.25) node [midway, above] {$60^\circ$};
	
  \coordinate (oo) at (0,1.25*\H);	
  \draw[stealth'-stealth'] ($(oo) + (0, 0.75)$) node[right] {\footnotesize $y$} -- (oo) -- ($(oo) + (0.75, 0)$) node[above] {\footnotesize $x$};


\end{tikzpicture}};
\end{tikzpicture}}
  \caption{Straight cantilever beam: Normalized vertical deflection versus beam slenderness for different in-plane mesh distortions (8 elements along the beam's length).}
  \label{fig:cantilever_distorted_lin}
\end{figure}
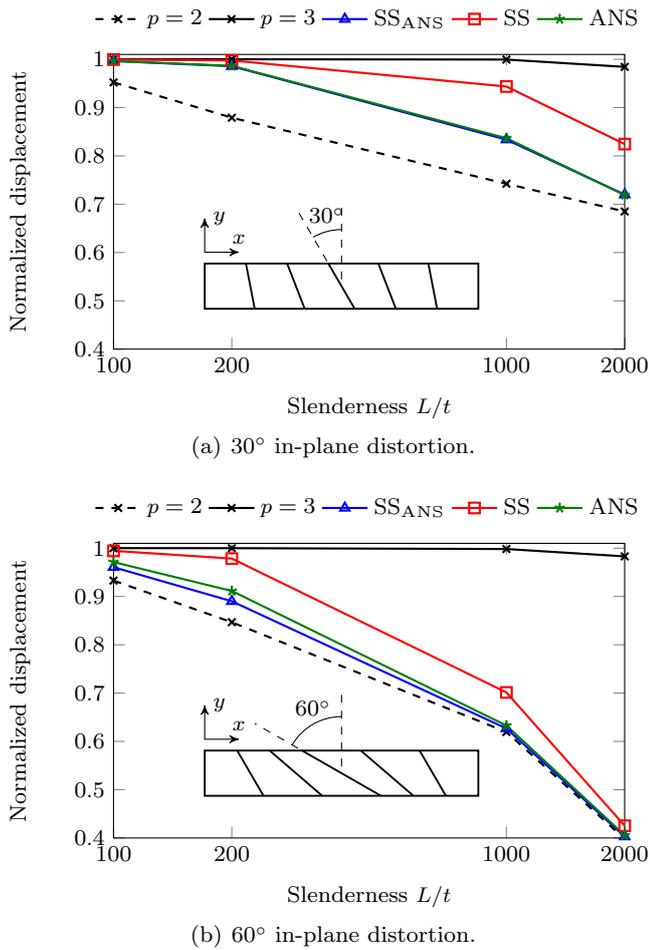

In order to test also the robustness of the proposed solid-shell elements, we report in Figure \ref{fig:cantilever_distorted_lin} the results obtained when different in-plane mesh distortions are considered.
These distortions are generated in such a way that the in-plane distortion angle is maximum at the center of the beam and linearly fades to zero at both ends (see Figure \ref{fig:cantilever_distorted_lin}), as proposed in \cite{Schwarze2011}.

As it can be seen in Figure \ref{fig:cantilever_distorted_lin_30}, for moderate angle distortions ($30^\circ$), the proposed Cartesian solid-shell element clearly outperforms the proposed ANS-inspired and ANS elements for large values of the slenderness, while for less slender beams, all elements present similar results.


For severe distortion angles ($60^\circ$), instead, Figure \ref{fig:cantilever_distorted_lin_60} shows that the Cartesian solid-shell element still provides better results than the other quadratic solid-shell formulations, but all of them attain quite inaccurate results for very slender beams ($L/t\geq1000$).

In all distorted cases the standard cubic element presents good results in the full considered range of slendernesses.

\subsection{Curved cantilever beam}
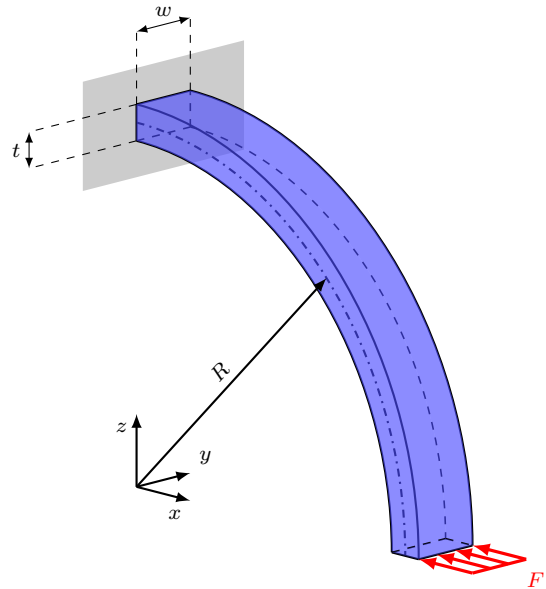
\begin{figure}
  \centering
  \tikzsetnextfilename{curved_geometry}\tdplotsetmaincoords{75}{45}
\begin{tikzpicture}[scale=0.5,tdplot_main_coords]
	
	\coordinate (O) at (0,0,0);
	\pgfmathsetmacro{\R}{10}
	\pgfmathsetmacro{\w}{2}
	\pgfmathsetmacro{\t}{1}
	\pgfmathsetmacro{\Rint}{\R-\t/2}
	\pgfmathsetmacro{\Rext}{\R+\t/2}
	
	\pgfmathsetmacro{\ax}{\Rext};
	\pgfmathsetmacro{\ay}{0};
	\pgfmathsetmacro{\bx}{\Rint};
	\pgfmathsetmacro{\by}{0};	
	
    \fill[fill=black!50!white,fill opacity=0.4] (0,-\w,-\t+\Rint) -- (0,2*\w,-\t+\Rint) -- (0,2*\w,2*\t+\Rint) -- (0,-\w,2*\t+\Rint) -- cycle;	

	\tdplotsetrotatedcoords{-90}{-90}{0}  
	\draw [line width=0.8pt, draw=black, tdplot_rotated_coords] (\Rint, 0, 0) arc (0:90:\Rint);
	\draw [line width=0.8pt, draw=black, tdplot_rotated_coords, dashdotted] (\R, 0, 0) arc (0:90:\R);
	\draw [line width=0.8pt, draw=black, tdplot_rotated_coords] (\Rext,0, 0) arc (0:90:\Rext);
	\draw [line width=0.6pt, draw=black,dashed, tdplot_rotated_coords] (\Rint, 0, \w) arc (0:90:\Rint); 
	\draw [line width=0.8pt, draw=black, tdplot_rotated_coords] (\Rext,0, \w) arc (0:90:\Rext);
	\draw [line width=0.8pt, draw=black, tdplot_rotated_coords] (\Rint, 0, 0) -- (\Rext, 0, 0);
	\draw [line width=0.6pt, draw=black, dashed, tdplot_rotated_coords] (\Rint, 0, \w) -- (\Rext, 0, \w);

	\draw [line width=0.8pt, draw=black] (0,0, \Rext) -- (0, \w, \Rext);
	\draw [line width=0.6pt, draw=black, dashed] (0,0, \Rint) -- (0, \w, \Rint);

	\tdplotsetrotatedcoords{-90}{-90}{0}
	\fill[fill=blue!70!white,fill opacity=0.4, tdplot_rotated_coords] (0, \Rext,0) arc (90:0:\Rext) --  (\ax, \ay, 0) -- (\bx, \by, 0) --  (\bx, \by, 0)  arc (0:90:\Rint) -- cycle;
	\fill[fill=blue!70!white,fill opacity=0.4, tdplot_rotated_coords, shift={(0,0,\w)}] (0, \Rext,0) arc (90:0:\Rext) --  (\ax, \ay, 0) -- (\bx, \by, 0) --  (\bx, \by, 0)  arc (0:90:\Rint) -- cycle;
	\fill[fill=blue!70!white,fill opacity=0.4, tdplot_rotated_coords] (0, \Rext,0) arc (90:0:\Rext) --  (\ax, \ay, 0) -- (\ax, \ay, \w) --  (\ax, \ay, \w)  arc (0:90:\Rext) -- cycle;
	\fill[fill=blue!70!white,fill opacity=0.4, tdplot_rotated_coords] (0, \Rint,0) arc (90:0:\Rint) --  (\bx, \by, 0) -- (\bx, \by, \w) --  (\bx, \by, \w)  arc (0:90:\Rint) -- cycle;

	\tdplotsetrotatedcoords{0}{0}{0}
	\draw [line width=0.8pt, draw=black, tdplot_rotated_coords] (\Rint, 0, 0) -- (\Rext, 0, 0);
	\draw [line width=0.6pt, draw=black, dashed, tdplot_rotated_coords] (\Rint, \w, 0) -- (\Rext, \w, 0);
	\draw [line width=0.6pt, draw=black, dashed, tdplot_rotated_coords] (\Rint, 0, 0) -- (\Rint, \w, 0);
	\draw [line width=0.8pt, draw=black, tdplot_rotated_coords] (\Rext, 0, 0) -- (\Rext, \w, 0);

	\fill[fill=blue!70!white,fill opacity=0.4, tdplot_rotated_coords] (\Rint, 0, 0) -- (\Rext, 0, 0) -- (\Rext, \w, 0) -- (\Rint, \w, 0) -- cycle;
	\fill[fill=blue!70!white,fill opacity=0.4] (0, 0, \Rint) -- (0, 0, \Rext) -- (0, \w, \Rext) -- (0, \w, \Rint) -- cycle;
	
    \draw[line width=0.8pt,->] (O) -- (2,0,0) node[anchor=north east]{$x$}; 
    \draw[line width=0.8pt,->] (O) -- (0,2,0) node[anchor=south west]{$y$};
    \draw[line width=0.8pt,->] (O) -- (0,0,2) node[anchor=north east]{$z$};	
    
	\begin{scope}[->, ultra thick]
    \draw[->,line width=1.2pt,red] (\Rext+2,0,0)--(\Rext,0,0);
    \draw[->,line width=1.2pt,red] (\Rext+2,\w/3,0)--(\Rext,\w/3,0);
    \draw[->,line width=1.2pt,red] (\Rext+2,2*\w/3,0)--(\Rext,2*\w/3,0);
    \draw[->,line width=1.2pt,red] (\Rext+2,\w,0)--(\Rext,\w,0);
  \end{scope}	
  \draw[-, line width=1.2pt,red] (\Rext+2,0,0)--(\Rext+2,\w,0);
    
    \draw[line width=0.4pt, draw=black, dashed] (0,0,\Rext) -- (0,0,\Rext+\w);
	\draw[line width=0.4pt, draw=black, dashed] (0, \w, \Rext) -- (0, \w, \Rext+\w);    
    \draw [<->,color=black] (0,0,\Rext+\w) -- (0, \w, \Rext+\w) node[black,midway,above] {$w$};
    \draw[line width=0.4pt, draw=black, dashed] (0,0,\Rint) -- (0,-2*\w, \Rint);
	\draw[line width=0.4pt, draw=black, dashed] (0,0,\Rext) -- (0,-2*\w, \Rext);    
    \draw [<->,color=black] (0,-2*\w, \Rint) -- (0, -2*\w, \Rext) node[black,sloped,midway,left, rotate=-90] {$t$};  
	\tdplotsetrotatedcoords{0}{-45}{0}   
	\draw [->, thick, tdplot_rotated_coords] (O) -- (\R, 0, 0) node [black,sloped,midway,above] {$R$}; 
    
    \draw [color=red] (\Rext+4,\w/2,0) node[anchor=east, rotate=0] {$F$};

\end{tikzpicture}
  \caption{Curved cantilever beam: Problem description.
	One face is clamped whereas a radial distributed load (with resultant $F$) is applied along the exterior edge
	of the opposite face. The radial beam deflection is measured at the interior edge of the free end of the beam.
  $E=1000$, $\nu=0$, $R=10$, and $w=1$. Different slendernesses $R/t$ are considered.}\label{fig:curved_geometry}
\end{figure}

The second example is a classical benchmark for membrane locking (see \cite{Echter2013}).
The problem is sketched in Figure \ref{fig:curved_geometry} and consists of a quarter of ring beam clamped at
one end and subjected to a radial distributed force
along the exterior edge of the opposite face, whose resultant is $F$.
The radius of the beam's middle fiber is $R=10$ and the width is $w=1$. As before, the Young's modulus and Poisson's ratio are $E=1000$ and $\nu=0$, respectively.

\begin{figure}
  \centering
  \tikzsetnextfilename{curved_beam_slenderness}\begin{tikzpicture}
 \begin{loglogaxis}[basic_paper_plots_style,
  xlabel={Slenderness $R/t$},
  ylabel={Normalized displacement},
  ymax=1.05,
  xmin=10, xmax=10000,
scale=0.45,
 legend columns=5,
 ]
 %
 %
 \addplot[solid 2]
 table [x=slenderness, y=S2]
 {tikz_figures/curved_beam_nel_10.res};
 \addlegendentry{\Stwo}
 \addplot[solid 3]
 table [x=slenderness, y=S3]
 {tikz_figures/curved_beam_nel_10.res};
 \addlegendentry{\Sthree}
 \addplot[ss ANS-like]
 table [x=slenderness, y=SS2-ANS-LIKE]
 {tikz_figures/curved_beam_nel_10.res};
 \addlegendentry{\SSANStwo}
 \addplot[ss]
 table [x=slenderness, y=SS2]
 {tikz_figures/curved_beam_nel_10.res};
 \addlegendentry{\SStwo}
 \addplot[ANS]
 table [x=slenderness, y=ANS2]
 {tikz_figures/curved_beam_nel_10.res};
 \addlegendentry{\ANS}
 \end{loglogaxis}
\end{tikzpicture}
  \caption{Curved cantilever beam: Normalized radial deflection for different slenderness values (10 elements along the beam's length).}\label{fig:curved_slen_lin}
\end{figure}
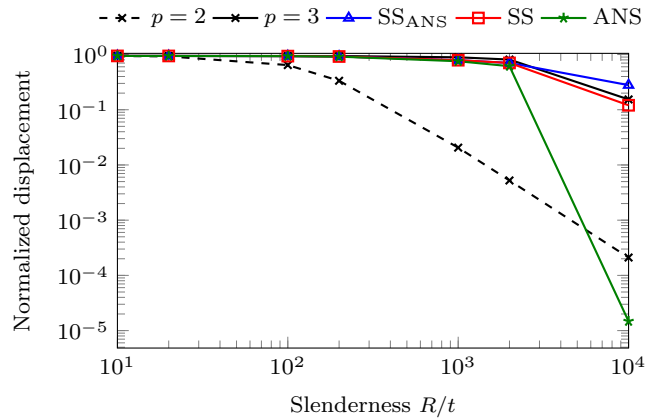

A mesh consisting of a single element for the beam cross section and 10 elements along the circumferential direction is considered.
The computed radial displacement at the beam tip is reported in Figure \ref{fig:curved_slen_lin} as a function of the slenderness $R/t$.
As it can be seen, all considered solid-shell elements significantly enhance the response of the standard quadratic formulation and all present very similar results that are practically locking-free, except for the case of extremely slender beams ($R/t=10^4$) for which the standard ANS formulation produces poor results.

\subsection{Shell obstacle course I: Scordelis-Lo roof}
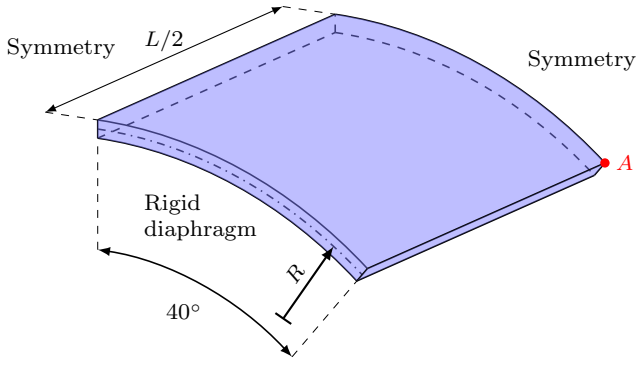
\begin{figure}
  \centering
  \tikzsetnextfilename{scordelis_loroof_geometry}
%
%
%
\tdplotsetmaincoords{75}{30}
\begin{tikzpicture}[scale=0.25,tdplot_main_coords]

\coordinate (O) at (0,0,0);
\pgfmathsetmacro{\R}{25};
\pgfmathsetmacro{\L}{25};
\pgfmathsetmacro{\t}{1};
\pgfmathsetmacro{\Rint}{\R-\t/2};
\pgfmathsetmacro{\Rext}{\R+\t/2};

\pgfmathsetmacro{\ax}{\Rext*cos(50)};
\pgfmathsetmacro{\ay}{\Rext*sin(50)};
\pgfmathsetmacro{\bx}{\Rint*cos(50)};
\pgfmathsetmacro{\by}{\Rint*sin(50)};

\tdplotsetrotatedcoords{-90}{-90}{0}  
\draw [line width=0.6pt, draw=black, tdplot_rotated_coords] (\Rint, 0, 0) arc (0:40:\Rint);
\draw [line width=0.6pt, draw=black, tdplot_rotated_coords, dashdotted] (\R, 0, 0) arc (0:40:\R);
\draw [line width=0.6pt, draw=black, tdplot_rotated_coords] (\Rext,0, 0) arc (0:40:\Rext);
\draw [line width=0.6pt, draw=black,dashed, tdplot_rotated_coords] (\Rint, 0, \L) arc (0:40:\Rint); 
\draw [line width=0.6pt, draw=black, tdplot_rotated_coords] (\Rext,0, \L) arc (0:40:\Rext);
\draw [line width=0.6pt, draw=black, tdplot_rotated_coords] (\Rint, 0, 0) -- (\Rext, 0, 0);
\draw [line width=0.6pt, draw=black, dashed, tdplot_rotated_coords] (\Rint, 0, \L) -- (\Rext, 0, \L);

\draw [line width=0.6pt, draw=black] (0,0, \Rext) -- (0, \L, \Rext);
\draw [line width=0.6pt, draw=black, dashed] (0,0, \Rint) -- (0, \L, \Rint);

\tdplotsetrotatedcoords{-90}{-90}{-50}
\fill[fill=blue!50!white,fill opacity=0.3, tdplot_rotated_coords] (0, \Rext,0) arc (90:50:\Rext) --  (\ax, \ay, 0) -- (\bx, \by, 0) --  (\bx, \by, 0)  arc (50:90:\Rint) -- cycle;
\fill[fill=blue!50!white,fill opacity=0.3, tdplot_rotated_coords, shift={(0,0,\L)}] (0, \Rext,0) arc (90:50:\Rext) --  (\ax, \ay, 0) -- (\bx, \by, 0) --  (\bx, \by, 0)  arc (50:90:\Rint) -- cycle;
\fill[fill=blue!50!white,fill opacity=0.3, tdplot_rotated_coords] (0, \Rext,0) arc (90:50:\Rext) --  (\ax, \ay, 0) -- (\ax, \ay, \L) --  (\ax, \ay, \L)  arc (50:90:\Rext) -- cycle;
\fill[fill=blue!50!white,fill opacity=0.3, tdplot_rotated_coords] (0, \Rint,0) arc (90:50:\Rint) --  (\bx, \by, 0) -- (\bx, \by, \L) --  (\bx, \by, \L)  arc (50:90:\Rint) -- cycle;

\tdplotsetrotatedcoords{0}{-50}{0}
\draw [line width=0.6pt, draw=black, tdplot_rotated_coords] (\Rint, 0, 0) -- (\Rext, 0, 0);
\draw [line width=0.6pt, draw=black, tdplot_rotated_coords] (\Rint, \L, 0) -- (\Rext, \L, 0);
\draw [line width=0.6pt, draw=black, tdplot_rotated_coords] (\Rint, 0, 0) -- (\Rint, \L, 0);
\draw [line width=0.6pt, draw=black, tdplot_rotated_coords] (\Rext, 0, 0) -- (\Rext, \L, 0);

\fill[fill=blue!50!white,fill opacity=0.3, tdplot_rotated_coords] (\Rint, 0, 0) -- (\Rext, 0, 0) -- (\Rext, \L, 0) -- (\Rint, \L, 0) -- cycle;
\fill[fill=blue!50!white,fill opacity=0.3] (0, 0, \Rint) -- (0, 0, \Rext) -- (0, \L, \Rext) -- (0, \L, \Rint) -- cycle;

\node[red, above=10pt, right=1pt, tdplot_rotated_coords] at (\Rext, \L, 0) {$A$};
\fill[fill=red, tdplot_rotated_coords] (\Rext, \L,0) circle (7.0pt);

\tdplotsetrotatedcoords{0}{-50}{0}   
\draw [line width=0.2pt, dashed, tdplot_rotated_coords] (3*\Rint/4, 0, 0) -- (\Rint, 0, 0); 
\tdplotsetrotatedcoords{0}{-90}{0}   
\draw [line width=0.2pt, dashed, tdplot_rotated_coords] (3*\Rint/4, 0, 0) -- (\Rint, 0, 0); 
\tdplotsetrotatedcoords{-90}{-90}{0}
\draw [<->, line width=0.6pt, tdplot_rotated_coords] (3*\Rint/4, 0, 0) arc (0:40:3*\Rint/4);
\node (O) [align=left] at (\R/4-1,0, 2*\R/3-1) {$40^\circ$};

\tdplotsetrotatedcoords{0}{-55}{0}   
\draw [|->, thick, tdplot_rotated_coords] (4*\Rint/5, 0, 0) -- (\R, 0, 0) node [black,midway,sloped,above] {$R$}; 

\draw [<->,color=black] (-\Rext/8,0, \Rext) -- (-\Rext/8, \L, \Rext) node[black,midway,above] {$L/2$};
\draw [line width=0.2pt, draw=black, dashed] (-\Rext/8,0, \Rext) -- (0,0, \Rext); 
\draw [line width=0.2pt, draw=black, dashed] (-\Rext/8,\L, \Rext) -- (0,\L, \Rext); 

\node (O) [align=left] at (\R-10, \L, \R) {Symmetry};
\node (O) [align=left] at (5, -\L/2, \R+8) {Symmetry};
\node (O) [align=left] at (\R/4, 0, \R-4) {Rigid\\ diaphragm};

\end{tikzpicture}
  \caption{Scordelis-Lo roof problem \cite{Belytschko1985}: Problem description. The cylindrical roof, supported by
	rigid diaphragms at both extrema is under the action of its self-weight $\rho g$. The displacement at mid-span point $A$ is measured.
	Due to symmetry conditions, only one quarter of the geometry is considered, setting the proper symmetry boundary conditions.
  $E=4.32\times 10^{8}$, $\nu=0$, $R=25$, $L=50$, $t=0.25$, and $\rho g = 360$.}\label{fig:scordelis_geometry}
\end{figure}

The first test of the shell obstacle course is the well-known Scordelis-Lo roof \cite{Scordelis1964}.
The roof, illustrated in Figure \ref{fig:scordelis_geometry}, has a cylindrical geometry and is supported by rigid diaphragms at both extremes.
The structure radius is $R=25$, its length is $L=50$, and the thickness is $t=0.25$.
The roof is subjected to its self-weight, whose
value is $\rho g = 360$, where $\rho$ is the density and $g$ is the gravity acceleration.
The elastic moduli are $E=4.32\cdot 10^{8}$ and $\nu=0$.
Due to symmetry conditions, only one quarter of the structure is modeled.  
The same number of elements are considered along each in-plane direction while only
one element through the thickness is used.

\begin{figure}
  \centering
  \tikzsetnextfilename{scordelis_convergence}\begin{tikzpicture}
 \begin{axis}[basic_paper_plots_style,
  xlabel={Number of control points per side},
  ylabel={Displacement},
  ymin=0.264, ymax=0.305,
  xmin=4, xmax=35,
scale=0.45
 ]
 %
 \addplot[reference] coordinates{
 (4, 0.3024)
 (35, 0.3024)}; 
 \addlegendentry{\reference\ ($0.3024$)}
%
%
%
%
 %
 %
 \addplot[solid 2]
 table [x=num_ctr_pts_side_2, y=S2]
 {tikz_figures/scordelis.res};
 \addlegendentry{\Stwo}
 \addplot[solid 3]
 table [x=num_ctr_pts_side_3, y=S3]
 {tikz_figures/scordelis.res};
 \addlegendentry{\Sthree}
 \addplot[ss ANS-like]
 table [x=num_ctr_pts_side_2, y=SS2-ANS-LIKE]
 {tikz_figures/scordelis.res};
 \addlegendentry{\SSANStwo}
 \addplot[ss]
 table [x=num_ctr_pts_side_2, y=SS2]
 {tikz_figures/scordelis.res};
 \addlegendentry{\SStwo}
 \addplot[ANS]
 table [x=num_ctr_pts_side_2, y=ANS2]
 {tikz_figures/scordelis.res};
 \addlegendentry{\ANS}
 \end{axis}
\end{tikzpicture}
  \caption{Scordelis-Lo roof problem: Vertical deflection at point $A$ versus number of in-plane control points per side.}\label{fig:scordelis_lin}
\end{figure}
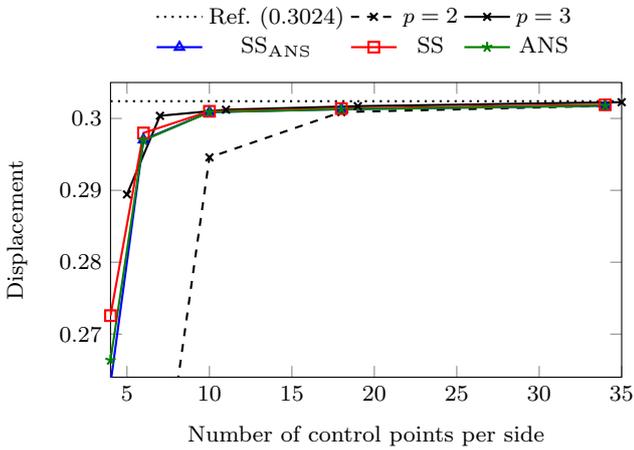

The reference displacement is the vertical deflection of  point $A$ in Figure \ref{fig:scordelis_geometry},
whose ``exact'' value is $0.3024$ (as reported in \cite{Belytschko1985}).
The results for the different elements considered, varying the number of in-plane control points, are reported in Figure \ref{fig:scordelis_lin}.
As it can be seen, the proposed solid-shell elements (as well as the ANS formulation) present results as good as those granted by the cubic discretization.

\subsection{Shell obstacle course II: Pinched hemispherical shell}
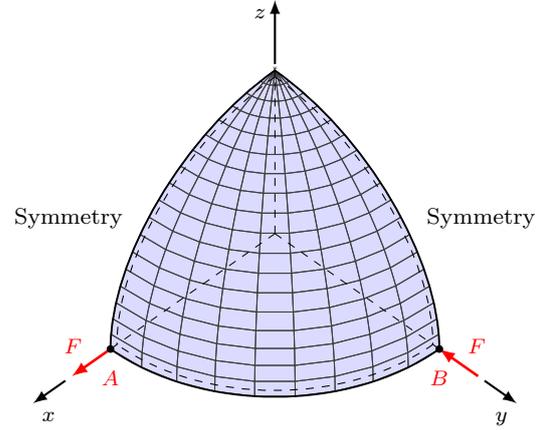
\begin{figure}
  \centering
  \tikzsetnextfilename{hemispherical_geometry}\tdplotsetmaincoords{135}{-45}
\begin{tikzpicture}[scale=0.12, tdplot_main_coords, fill opacity=.2,draw opacity=1]

\coordinate (O) at (0,0,0);
\pgfmathsetmacro{\R}{25};
\pgfmathsetmacro{\t}{1};
\pgfmathsetmacro{\Rint}{\R-\t/2};
\pgfmathsetmacro{\Rext}{\R+\t/2};

\tdplotsetpolarplotrange{0}{90}{0}{90};
\tdplotsphericalsurfaceplot[]{72}{36}{\Rext}{black!75!white}{blue!70!white}{}{}{};

\draw[->, line width=0.8pt, draw=black, fill = black, opacity=1] (\Rext+7,0,0) -- (\Rext+12,0,0) node[anchor=north east]{$y$}; 
\draw[->, line width=0.8pt, draw=black, fill = black, opacity=1] (0, \Rext+7,0) -- (0,\Rext+12,0) node[anchor=north west]{$x$};
\draw[->, line width=0.8pt, draw=black, fill = black, opacity=1] (0,0,\Rext+1) -- (0,0,\Rext+11) node[anchor=north east]{$z$};

\begin{scope}[->, ultra thick]
\draw[->,line width=1.0pt,red, opacity=1] (\Rext+6, 0, 0)--(\Rext,0,0) node[red,midway,above right] {$F$};
\draw[->,line width=1.0pt,red, opacity=1] (0, \Rext,0) -- (0, \Rext+6, 0) node[red,midway,above left] {$F$};
\end{scope}	

\node[red, opacity=1, below=5pt] at (0, \Rext,0) {$A$};
\fill[black, opacity=1] (0, \Rext,0) circle (12.0pt);

\node[red, opacity=1, below=5pt] at (\Rext, 0, 0) {$B$};
\fill[black, opacity=1] (\Rext, 0, 0) circle (12.0pt);

\draw [line width=0.4pt, draw=black,dashed] (\Rint, 0, 0) arc (0:90:\Rint);
\draw [line width=0.7pt, draw=black] (\Rext, 0, 0) arc (0:90:\Rext);
\draw [line width=0.4pt, draw=black,dashed] (\Rint, 0, 0) -- (\Rext, 0, 0);
\draw [line width=0.4pt, draw=black,dashed] (0, \Rint, 0) -- (0, \Rext, 0);  
\draw [line width=0.4pt, draw=black,dashed] (0,0, \Rint) -- (0, 0, \Rext);

\foreach \theta in {0,-90}  
{
\tdplotsetrotatedcoords{\theta}{-90}{0}   
\draw [line width=0.4pt, draw=black,dashed, tdplot_rotated_coords] (\Rint, 0, 0) arc (0:90:\Rint); 
\draw [line width=0.7pt, draw=black, tdplot_rotated_coords] (\Rext, 0, 0) arc (0:90:\Rext); 
}

\draw [line width=0.2pt, draw=black, dashed] (O) -- (\Rint, 0, 0);
\draw [line width=0.2pt, draw=black, dashed] (O) -- (0, \Rint, 0);
\draw [line width=0.2pt, draw=black, dashed] (O) -- (0, 0, \Rint);

\node (O) [align=left, black, opacity=1] at (\R+7, 0, \R) {Symmetry};
\node (O) [align=left, black, opacity=1] at (0, \R+7, \R) {Symmetry};

\end{tikzpicture}
  \caption{Pinched hemispherical shell \cite{Belytschko1985}: Problem description. The hemispherical structure, fixed at the top point, is subjected
	to the action of diametrically opposed forces (applied at points $A$ and $B$ in the picture).
	The equator (bottom) section can move freely. The radial deflection at point $A$ is measured.
	Due to symmetry conditions, only one quarter of the geometry is considered, setting the proper symmetry boundary conditions.
  $E=6.825\times10^{7}$, $\nu=0.3$, $R=10$, $t=0.04$ and $F=1$.}\label{fig:hemispherical_geometry}
\end{figure}

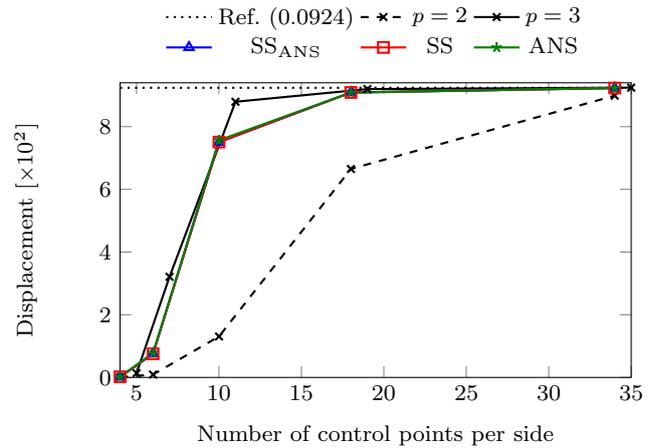
\begin{figure}
  \centering
  \tikzsetnextfilename{hemispherical_convergence}\begin{tikzpicture}
 \begin{axis}[basic_paper_plots_style,
  xlabel={Number of control points per side},
  ylabel={Displacement [$\times 10^2$]},
  ymin=0, ymax=9.4,
  xmin=4, xmax=35,
scale=0.45
 ]
 \addplot[reference] coordinates{
 (4, 9.24)
 (35, 9.24)}; 
 \addlegendentry{\reference\ ($0.0924$)}
 %
 %
 \addplot[solid 2]
 table [x=num_ctr_pts_side_2, y expr=\thisrow{S2}*1.0e2]
 {tikz_figures/hemispherical.res};
 \addlegendentry{\Stwo}
 \addplot[solid 3]
 table [x=num_ctr_pts_side_3, y expr=\thisrow{S3}*1.0e2]
 {tikz_figures/hemispherical.res};
 \addlegendentry{\Sthree}
 \addplot[ss ANS-like]
 table [x=num_ctr_pts_side_2, y expr=\thisrow{SS2-ANS-LIKE}*1.0e2]
 {tikz_figures/hemispherical.res};
 \addlegendentry{\SSANStwo}
 \addplot[ss]
 table [x=num_ctr_pts_side_2, y expr=\thisrow{SS2}*1.0e2]
 {tikz_figures/hemispherical.res};
 \addlegendentry{\SStwo}
 %
 %
 %
 \addplot[ANS]
 table [x=num_ctr_pts_side_2, y expr=\thisrow{ANS2}*1.0e2]
 {tikz_figures/hemispherical.res};
 \addlegendentry{\ANS}
 \end{axis}
\end{tikzpicture}
  \caption{Pinched hemispherical shell \cite{Belytschko1985}: Radial deflection at point A versus  the number of control points along each in-plane direction.}\label{fig:hemispherical_lin}
\end{figure}

The second test of the shell obstacle course consists in a hemispherical structure pinched by two couples of opposite concentrated forces on diametrically opposed points of the equator section.
Due to the problem symmetry only one quarter of the structure is modeled, setting suitable symmetry boundary conditions, as illustrated in Figure \ref{fig:hemispherical_geometry}.
Moreover, the top point of the hemisphere is fixed, while the equator section of the hemisphere can move freely.
The structure has a thickness $t=0.04$ and the middle surface of the hemisphere has a radius $R=10$. 
The applied forces have magnitude $F=1$ and the material properties are $E=6.825\times10^{7}$ and $\nu=0.3$.
In Figure \ref{fig:hemispherical_lin} we plot the numerically computed radial deflection at point $A$  (see Figure \ref{fig:hemispherical_geometry}) versus the number of control points along each in-plane direction and compare those results with the reference solution $0.0924$ reported in \cite{Belytschko1985}.
A similar behavior to that observed in the Scordelis-Lo roof case is obtained.

\subsection{Shell obstacle course III: Pinched cylinder}
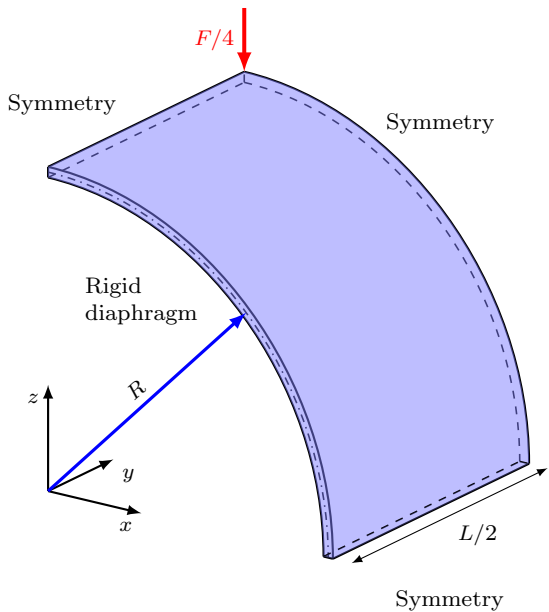
\begin{figure}
  \centering
  \tikzsetnextfilename{pinched_geometry}
%
%
%
\tdplotsetmaincoords{70}{35}
\begin{tikzpicture}[scale=0.15,tdplot_main_coords]
	
	\coordinate (O) at (0,0,0);
	\pgfmathsetmacro{\R}{30}
	\pgfmathsetmacro{\L}{30}
	\pgfmathsetmacro{\t}{1.0}
	\pgfmathsetmacro{\Rint}{\R-\t/2}
	\pgfmathsetmacro{\Rext}{\R+\t/2}

	\pgfmathsetmacro{\ax}{\Rext};
	\pgfmathsetmacro{\ay}{0};
	\pgfmathsetmacro{\bx}{\Rint};
	\pgfmathsetmacro{\by}{0};	
	
	\tdplotsetrotatedcoords{-90}{-90}{0}  
	\draw [line width=0.8pt, draw=black, tdplot_rotated_coords] (\Rint, 0, 0) arc (0:90:\Rint);
	\draw [line width=0.6pt, draw=black, tdplot_rotated_coords, dashdotted] (\R, 0, 0) arc (0:90:\R);
	\draw [line width=0.8pt, draw=black, tdplot_rotated_coords] (\Rext,0, 0) arc (0:90:\Rext);
	\draw [line width=0.6pt, draw=black,dashed, tdplot_rotated_coords] (\Rint, 0, \L) arc (0:90:\Rint); 
	\draw [line width=0.8pt, draw=black, tdplot_rotated_coords] (\Rext,0, \L) arc (0:90:\Rext);
	\draw [line width=0.8pt, draw=black, tdplot_rotated_coords] (\Rint, 0, 0) -- (\Rext, 0, 0);
	\draw [line width=0.6pt, draw=black, dashed, tdplot_rotated_coords] (\Rint, 0, \L) -- (\Rext, 0, \L);

	\draw [line width=0.8pt, draw=black] (0,0, \Rext) -- (0, \L, \Rext);
	\draw [line width=0.6pt, draw=black, dashed] (0,0, \Rint) -- (0, \L, \Rint);

	\tdplotsetrotatedcoords{-90}{-90}{0}
	\fill[fill=blue!50!white,fill opacity=0.3, tdplot_rotated_coords] (0, \Rext,0) arc (90:0:\Rext) --  (\ax, \ay, 0) -- (\bx, \by, 0) --  (\bx, \by, 0)  arc (0:90:\Rint) -- cycle;
	\fill[fill=blue!50!white,fill opacity=0.3, tdplot_rotated_coords, shift={(0,0,\L)}] (0, \Rext,0) arc (90:0:\Rext) --  (\ax, \ay, 0) -- (\bx, \by, 0) --  (\bx, \by, 0)  arc (0:90:\Rint) -- cycle;
	\fill[fill=blue!50!white,fill opacity=0.3, tdplot_rotated_coords] (0, \Rext,0) arc (90:0:\Rext) --  (\ax, \ay, 0) -- (\ax, \ay, \L) --  (\ax, \ay, \L)  arc (0:90:\Rext) -- cycle;
	\fill[fill=blue!50!white,fill opacity=0.3, tdplot_rotated_coords] (0, \Rint,0) arc (90:0:\Rint) --  (\bx, \by, 0) -- (\bx, \by, \L) --  (\bx, \by, \L)  arc (0:90:\Rint) -- cycle;

	\tdplotsetrotatedcoords{0}{0}{0}
	\draw [line width=0.8pt, draw=black, tdplot_rotated_coords] (\Rint, 0, 0) -- (\Rext, 0, 0);
	\draw [line width=0.6pt, draw=black, dashed, tdplot_rotated_coords] (\Rint, \L, 0) -- (\Rext, \L, 0);
	\draw [line width=0.6pt, draw=black, dashed, tdplot_rotated_coords] (\Rint, 0, 0) -- (\Rint, \L, 0);
	\draw [line width=0.8pt, draw=black, tdplot_rotated_coords] (\Rext, 0, 0) -- (\Rext, \L, 0);

	\fill[fill=blue!50!white,fill opacity=0.3, tdplot_rotated_coords] (\Rint, 0, 0) -- (\Rext, 0, 0) -- (\Rext, \L, 0) -- (\Rint, \L, 0) -- cycle;
	\fill[fill=blue!50!white,fill opacity=0.3] (0, 0, \Rint) -- (0, 0, \Rext) -- (0, \L, \Rext) -- (0, \L, \Rint) -- cycle;	
	
    \draw[line width=0.8pt,->] (O) -- (10,0,0) node[anchor=north east]{$x$}; 
    \draw[line width=0.8pt,->] (O) -- (0,10,0) node[anchor=north west]{$y$};
    \draw[line width=0.8pt,->] (O) -- (0,0,10) node[anchor=north east]{$z$};	
    
	\begin{scope}[->, ultra thick]
    \draw[->,line width=1.5pt,red] (0,\L,\Rext+6)--(0,\L,\Rext) node[red,midway,left] {$F/4$};
    \end{scope}	
    
	\draw [<->,color=black] (\Rext+2,0,0) -- (\Rext+2, \L, 0) node[black,midway,below right] {$L/2$};  
	\tdplotsetrotatedcoords{0}{-45}{0}   
	\draw [->,very thick, blue, tdplot_rotated_coords] (O) -- (\R, 0, 0) node [black,sloped,midway,above] {$R$}; 
	
	\node (O) [align=left] at (\R-9, \L, \R) {Symmetry};
	\node (O) [align=left] at (-9, \L/2, \R) {Symmetry};
	\node (O) [align=left] at (\Rext+2, \L/2, -8) {Symmetry};

	\node (O) [align=left] at (\R/3, 0, \R-10) {Rigid\\ diaphragm};
	
\end{tikzpicture}
%
  \caption{Pinched cylinder \cite{Belytschko1985}: Problem description.
  The cylindrical structure is supported by two rigid diaphragms
	at both ends while a pair of diametrically opposed concentrated forces
	are applied at the center of the cylinder. The radial deflection of the point
	where  loads are applied is measured.
	Due to symmetry conditions, only one eighth of the geometry is considered, setting the proper symmetry boundary conditions.
  $E=3\times 10^{6}$, $\nu=0.3$, $R=300$, $L=600$, $t=3$, and $F=1$.}\label{fig:pinched_geometry}
\end{figure}

The last test case of the shell obstacle course, is the so-called pinched cylinder.
The problem consists of a cylinder with rigid end diaphragms subjected to a pair of concentrated forces; due to the problem symmetry only one eight of the problem is studied, setting the proper symmetry boundary conditions, as shown in Figure \ref{fig:pinched_geometry}. 
The cylinder middle surface has radius $R=300$, while its length is $L=600$ and its thickness $t=3$.
The two opposite concentrated forces have  value $F=1$, and the material properties are $E=3\times 10^{6}$ and $\nu=0.3$.

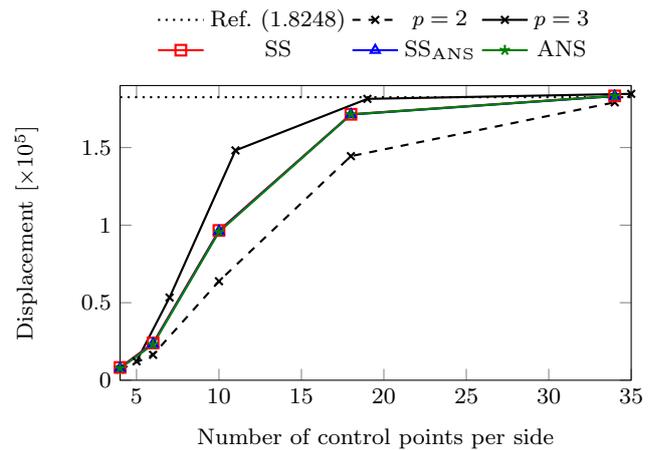
\begin{figure}
  \centering
  \tikzsetnextfilename{cylinder_convergence}\begin{tikzpicture}
 \begin{axis}[basic_paper_plots_style,
  xlabel={Number of control points per side},
  ylabel={Displacement [$\times 10^5$]},
  xmin=4, xmax=35,
 ymin=0, ymax=1.9,
 legend columns=3,
 scale=0.45
 ]
 \addplot[reference] coordinates{
 (4,  1.8248)
 (35, 1.8248)}; 
 \addlegendentry{\reference\ ($1.8248$)}
 %
 %
 \addplot[solid 2]
 table [x=num_ctr_pts_side_2, y expr=\thisrow{S2}*1.0e5]
 {tikz_figures/cylinder.res};
 \addlegendentry{\Stwo}
 \addplot[solid 3]
 table [x=num_ctr_pts_side_3, y expr=\thisrow{S3}*1.0e5]
 {tikz_figures/cylinder.res};
 \addlegendentry{\Sthree}
 \addplot[ss]
 table [x=num_ctr_pts_side_2,  y expr=\thisrow{SS2}*1.0e5]
 {tikz_figures/cylinder.res};
 \addlegendentry{\SStwo}
 \addplot[ss ANS-like]
 table [x=num_ctr_pts_side_2,  y expr=\thisrow{SS2-ANS-LIKE}*1.0e5]
 {tikz_figures/cylinder.res};
 \addlegendentry{\SSANStwo}
 \addplot[ANS]
 table [x=num_ctr_pts_side_2,   y expr=\thisrow{ANS2}*1.0e5]
 {tikz_figures/cylinder.res};
 \addlegendentry{\ANS}
 \end{axis}
\end{tikzpicture}
  \caption{Pinched cylinder: Radial deflection versus the number control points along each in-plane direction.}\label{fig:pinched_lin}
\end{figure}

The numerically computed radial deflections (measured at the point where  loads are applied) are plotted in Figure \ref{fig:pinched_lin} versus the number of control points along each in-plane direction and compared with the reference solution $1.8248\times 10^{-5}$ reported in \cite{Belytschko1985}.
As in the previous test cases, all  quadratic solid-shell elements present similar behaviors, representing a very significant improvement over the standard quadratic element. In this case,  we have to note however that the standard cubic element shows a slightly superior performance.

\section{Conclusions}
In this paper, we have presented a new approach to alleviate geometrical locking effects in solid shells.
The approach is based on local projections of strains onto coarser polynomial spaces. We have explored two different formulations based on this method.
The first one is inspired by the ANS method and uses different projection spaces for the different strain components, while in the second formulation, the same projection is used for all strain components.
Both formulations have shown very good performance in all numerical tests with the same level of accuracy as the ANS formulation presented in \cite{Caseiro2014}.
The advantages of the proposed formulations are their simplicity and numerical efficiency, requiring much fewer function evaluations at the element level than the ANS method.
Comparing the two presented formulations, the second one is even simpler and more efficient than the first one.
This formulation requires only a standard Cartesian-based element formulation and can be integrated into existing solid implementations very easily.
In this paper, we restricted to linear elastic problems.
The extension to nonlinear mechanics is planned as future work.

\bibliographystyle{spmpsci}      
\bibliography{references}

\begin{acknowledgements}
Pablo Antol\'in gratefully acknowledges the support of the European Research Council, through the ERC AdG n.\ 694515 - CHANGE grant,
and Alessandro Reali has been also partially supported by the MIUR-PRIN project XFAST-SIMS (no.\ 20173C478N).
\end{acknowledgements}

\appendix
\section{Local $L^2$ projection operators} \label{sec:proj_coeffs}
In this Appendix we provide closed-form expressions for the projection operators $\mathbb{P}^{(i,j)}$ introduced in Section \ref{sec:ANS-inspired}, for degrees $p=1$ and $p=2$.
For both degrees, $p+1$ Gauss-Legendre quadrature points along each direction are considered, being the points ordered is a lexicographical manner: i.e., the first parametric direction runs faster than the second, and the second faster than the third one.

The projection operator matrices present the following block-diagonal structure:
\begin{align}
  \mathbb{P}^{(i,j)} &= \underbrace{\begin{bmatrix}
	\mathbb{S}^{(i,j)} & \mathbf{0}     & \dots  &  \mathbf{0}     \\
	\mathbf{0}     & \mathbb{S}^{(i,j)} & \dots  &  \mathbf{0}     \\
	\vdots          & \vdots        & \ddots &  \vdots         \\
	\mathbf{0}     & \mathbf{0}     & \dots  &  \mathbb{S}^{(i,j)}
  \end{bmatrix}}_{(p+1)\text{ submatrices}}\,,
\end{align}
where the submatrices $\mathbb{S}^{(i,j)}\in\mathbb{R}^{(p+1)^2\times(p+1)^2}$ for degree $p=1$ are:
\begin{subequations}
\begin{align}
  \mathbb{S}^{(1,1)} &= \frac{1}{2} \begin{bmatrix}
1 & 1 & 0 & 0 \\
1 & 1 & 0 & 0 \\
0 & 0 & 1 & 1 \\
0 & 0 & 1 & 1
\end{bmatrix}\,,\\
\mathbb{S}^{(2,2)} &= \frac{1}{2} \begin{bmatrix}
1 & 0 & 1 & 0 \\
0 & 1 & 0 & 1 \\
1 & 0 & 1 & 0 \\
0 & 1 & 0 & 1
  \end{bmatrix}\,,\\
\mathbb{S}^{(1,2)} &= \frac{1}{4} \begin{bmatrix}
1 & 1 & 1 & 1 \\
1 & 1 & 1 & 1 \\
1 & 1 & 1 & 1 \\
1 & 1 & 1 & 1
  \end{bmatrix}\,,
\end{align}
\end{subequations}
and for degree $p=2$:
\begin{subequations}
\begin{align}
  \mathbb{S}^{(1,1)} &= \frac{1}{18} \begin{bmatrix}
    14 & 8 & -4 & 0 & 0 & 0 & 0 & 0 & 0 \\
    5 & 8 & 5 & 0 & 0 & 0 & 0 & 0 & 0 \\
    -4 & 8 & 14 & 0 & 0 & 0 & 0 & 0 & 0 \\
    0 & 0 & 0 & 14 & 8 & -4 & 0 & 0 & 0 \\
    0 & 0 & 0 & 5 & 8 & 5 & 0 & 0 & 0 \\
    0 & 0 & 0 & -4 & 8 & 14 & 0 & 0 & 0 \\
    0 & 0 & 0 & 0 & 0 & 0 & 14 & 8 & -4 \\
    0 & 0 & 0 & 0 & 0 & 0 & 5 & 8 & 5 \\
    0 & 0 & 0 & 0 & 0 & 0 & -4 & 8 & 14
		\end{bmatrix}\,,\\
  \mathbb{S}^{(2,2)} &= \frac{1}{18} \begin{bmatrix}
14 & 0 & 0 & 8 & 0 & 0 & -4 & 0 & 0 \\
0 & 14 & 0 & 0 & 8 & 0 & 0 & -4 & 0 \\
0 & 0 & 14 & 0 & 0 & 8 & 0 & 0 & -4 \\
5 & 0 & 0 & 8 & 0 & 0 & 5 & 0 & 0 \\
0 & 5 & 0 & 0 & 8 & 0 & 0 & 5 & 0 \\
0 & 0 & 5 & 0 & 0 & 8 & 0 & 0 & 5 \\
-4 & 0 & 0 & 8 & 0 & 0 & 14 & 0 & 0 \\
0 & -4 & 0 & 0 & 8 & 0 & 0 & 14 & 0 \\
0 & 0 & -4 & 0 & 0 & 8 & 0 & 0 & 14
  \end{bmatrix}\,,\\
\mathbb{S}^{(1,2)} &= \scriptsize\frac{1}{324} \begin{bmatrix}
196 & 112 & -56 & 112 & 64 & -32 & -56 & -32 & 16 \\
70 & 112 & 70 & 40 & 64 & 40 & -20 & -32 & -20 \\
-56 & 112 & 196 & -32 & 64 & 112 & 16 & -32 & -56 \\
70 & 40 & -20 & 112 & 64 & -32 & 70 & 40 & -20 \\
25 & 40 & 25 & 40 & 64 & 40 & 25 & 40 & 25 \\
-20 & 40 & 70 & -32 & 64 & 112 & -20 & 40 & 70 \\
-56 & -32 & 16 & 112 & 64 & -32 & 196 & 112 & -56 \\
-20 & -32 & -20 & 40 & 64 & 40 & 70 & 112 & 70 \\
16 & -32 & -56 & -32 & 64 & 112 & -56 & 112 & 196
  \end{bmatrix}\,.
\end{align}
\end{subequations}
\end{document}